\documentclass{amsart}
\usepackage{graphicx}
\usepackage{amssymb} \usepackage{epstopdf} 
\usepackage{nicefrac}
\usepackage{enumerate}
\usepackage{hyperref}
\usepackage{float}
\usepackage{color}
\usepackage{pst-all}
\usepackage{tabularx}
\usepackage{soul}

\newtheorem{thm}{THEOREM}[section]

\newtheorem{cor}[thm]{COROLLARY}

\newtheorem{defn}[thm]{DEFINITION}

\newtheorem{lemma}[thm]{LEMMA}
\newtheorem{prob}[thm]{PROBLEM}
\newtheorem{prop}[thm]{PROPOSITION}

\newtheorem{problem}[thm]{PROBLEM}
\newtheorem{remark}[thm]{REMARK}

\newcommand{\ds}{\displaystyle}

\newcommand{\cD}{{\mathcal D}}
\newcommand{\cG}{{\mathcal G}}
 
\newcommand{\cH}{{\mathcal H}}

\newcommand{\cN}{{\mathcal N}}
\newcommand{\cO}{{\mathcal O}}

\newcommand{\cR}{{\mathcal R}}

\newcommand{\CO}{{\rm CO}} 

\newcommand{\cS}{{\mathcal S}}
\newcommand{\cU}{{\mathcal U}}

\newcommand{\diam}{{\rm diam}} 
\newcommand{\dX}{d_{\fX}} 
\newcommand{\e}{{\varepsilon}} 

\newcommand{\fD}{{\mathfrak{D}}}
\newcommand{\fG}{{\mathfrak{G}}}
\newcommand{\fH}{{\mathfrak{H}}}
\newcommand{\fN}{{\mathfrak{N}}}
\newcommand{\Fix}{{\rm Fix}}

\newcommand{\fX}{{\mathfrak{X}}}

\newcommand{\Iso}{{\rm Iso}}

\newcommand{\mZ}{{\mathbb Z}}

\newcommand{\vp}{{\varphi}}

 \newcommand{\Homeo}{{\rm Homeo}} 
 \newcommand{\Perm}{{\rm Perm}}

\newcommand{\Aut}{{\rm Aut}}

\newcommand{\whPhi}{\widehat{\Phi}}

\newcommand{\whId}{\widehat{\rm Id}}
\newcommand{\whe}{\widehat{e}}
\newcommand{\whg}{\widehat{g}}

\newcommand{\whx}{\widehat{x}}

\newcommand{\whA}{\widehat{A}}

\newcommand{\whN}{\widehat{N}}

\newcommand{\whU}{\widehat{U}}

\newcommand{\whGamma}{{\widehat{\Gamma}}}

\newcommand{\whtau}{{\widehat{\tau}}}
\newcommand{\whvp}{{\widehat{\varphi}}}

\newcommand{\G}{\Gamma}

\input xy
\xyoption {all}

\begin{document}

\title{Cantor dynamics of renormalizable groups}

\author{Steven Hurder}
\address{Steven Hurder, Department of Mathematics, University of Illinois at Chicago, 322 SEO (m/c 249), 851 S. Morgan Street, Chicago, IL 60607-7045}
\email{hurder@uic.edu}

\author{Olga Lukina}
\address{Olga Lukina, Faculty of Mathematics, University of Vienna, Oskar-Morgenstern-Platz 1, 1090 Vienna, Austria}
\email{olga.lukina@univie.ac.at}
  
\author{Wouter Van Limbeek}
\address{Wouter Van Limbeek, Department of Mathematics, University of Illinois at Chicago, 322 SEO (m/c 249), 851 S. Morgan Street, Chicago, IL 60607-7045}
\email{wouter@uic.edu}

\thanks{Version date: February 4, 2020; revision October 28, 2020}

\thanks{2010 {\it Mathematics Subject Classification}. Primary:  20F18, 37B50, 37C85}
 
\begin{abstract}
A group $\G$ is said to be finitely non-co-Hopfian, or renormalizable, if there exists   a self-embedding $\vp \colon \G \to \G$ whose image is a proper subgroup of  finite index. Such a proper self-embedding is called a renormalization for $\G$. In this work, we associate  a    dynamical system to a renormalization $\vp$ of $\G$. The discriminant invariant $\cD_{\vp}$ of the associated Cantor dynamical system is a profinite group which  is a measure of the asymmetries of the dynamical system.  If $\cD_{\vp}$  is a finite group for some renormalization, we show that  $\G/C_{\vp}$ is virtually nilpotent, where $C_{\vp}$ is the kernel of the action map.  We    introduce the notion of a (virtually) renormalizable Cantor action,   and show that  the action associated to a renormalizable group is   virtually renormalizable.  We study the  properties of  virtually renormalizable    Cantor actions, and show that virtual renormalizability   is an invariant of continuous orbit equivalence. Moreover, the discriminant invariant of a renormalizable Cantor action is   an invariant of continuous orbit equivalence. Finally, the notion of a renormalizable Cantor action is related to the notion of a self-replicating group of automorphisms of a rooted tree. 
 \end{abstract}

  \thanks{Keywords: non-co-Hopfian groups, minimal   Cantor actions, odometer actions, renormalization}

  \maketitle

\section{Introduction}\label{sec-intro}
  
A countable group  $\G$   is \emph{co-Hopfian}   if every monomorphism $\vp \colon \G \to \G$ is an isomorphism   \cite{Baer1944}, and    is said to be \emph{non-co-Hopfian} otherwise. If there exists   a self-embedding $\vp$ whose image is a proper subgroup of  finite index, then $\G$ is said to be \emph{finitely non-co-Hopfian} \cite{WvL2018b}.  A proper self-embedding $\vp \colon \G \to \G$ with finite index is called a \emph{renormalization} of $\G$, in analogy with the case for    $\G = \mZ^n$. If $\G$ admits a renormalization, then  it is  said to be  \emph{renormalizable}.

  The free abelian group $\mZ^n$ is   renormalizable, as are   many  finitely generated nilpotent groups. There are  also many examples of renormalizable groups  which are not nilpotent, as   described for example in \cite{DP2003,ER2005,GW1992,GW1994,GLiW1994,NekkyPete2011,OP1998,WvL2018b}.  On the other hand, the free group $\mZ^{\star n} = \mZ \star  \cdots \star \mZ$ for $n \geq 2$   is non-co-Hopfian, but is not renormalizable. The classification of   non-co-Hopfian groups in general appears to be a difficult problem. 
  
  There is a related concept of a \emph{scale-invariant}   group,  introduced by Benjamini (see \cite[Section~9.2]{Sapir2007}). A \emph{scale} for $\G$ is a descending chain   of finite index subgroups $\cS = \{\G_{\ell} \mid \ell \geq 1\}$ whose intersection is a finite group, and such that for each $\ell$, there exists an isomorphism $\phi_{\ell} \colon \G \to \G_{\ell}$. Benjamini asked if a scale-invariant group must be virtually nilpotent? Nekrashevych and Pete \cite[Theorem~1.1]{NekkyPete2011} gave examples   of scale-invariant groups which are not virtually nilpotent. In the same work, the authors defined the notion of a \emph{strongly scale-invariant} group, as a renormalizable group $\G$ such that the collection of subgroups $\{\G_{\ell} = \vp^{\ell}(\G)  \mid \ell \geq 0\}$ is a scale for $\G$. Then \cite[Question~1.1]{NekkyPete2011} asks if a strongly scale-invariant group must be virtually nilpotent?

 In this work, we obtain   a partial answer to this question. The key idea is to study the properties of an infinite 
  profinite group $\whGamma_{\vp}$ naturally associated to a renormalization $\vp$. The group $\whGamma_{\vp}$ is a   proper quotient of the full profinite completion of $\G$. The key observation is given   in Proposition~\ref{prop-embedding}, which shows  that $\vp \colon \G \to \G$ induces   an \emph{open} embedding $\whvp \colon \whGamma_{\vp} \to \whGamma_{\vp}$.     Then Theorem~\ref{thm-contractionprinciple} states that if both the intersections    
  $\ds \cap_{\ell > 0} \ \vp^{\ell}(\G)$ and  $\ds \cap_{\ell > 0} \ \whvp^{\ell}(\whGamma_{\vp})$ are      finite groups, then $\G$ is  virtually nilpotent.     In other words, we answer in the affirmative the question of Nekrashevych and Pete above under a stronger assumption, that both $\G$ and the infinite profinite group $\whGamma_\vp$ admit a scale.  
  
 More precisely, given a renormalization $\vp \colon \G \to \G$,    let $C(\cG_{\vp})$ be the largest normal subgroup  of   $K(\cG_{\vp}) = \bigcap_{\ell > 0} \vp^{\ell}(\G)$.       In Definition~\ref{def-discriminant} we define a closed subgroup       $\cD_{\vp} \subset \whGamma_{\vp}$ which is naturally associated to the renormalization $\vp$. In Theorem \ref{thm-main2}, we prove that $\cD_{\vp}=\bigcap_{\ell>0} \whvp^\ell(\whGamma_{\vp})$. In particular, $\vp$ induces a scale on the profinite group $\whGamma_{\vp}$ if and only if $\cD_{\vp}$ is finite. We have:
          
\begin{thm}\label{thm-main4}
 Let $\G$ be a finitely generated group, and let $\vp \colon \G \to \G$ be a renormalization of $\G$.
 \begin{enumerate}
\item If  $\cD_{\vp}$   is the trivial group,    then    $\G/C(\cG_{\vp})$ is nilpotent. 
\item If $\cD_{\vp}$   is a finite group,    then    $\G/C(\cG_{\vp})$ is virtually nilpotent.
\end{enumerate}
    \end{thm}
The assumption that $\G$ is finitely generated is essential, as shown by the examples in Section~\ref{subsec-infinite}.

Our approach to the study of  renormalizable groups  is based on the study of the Cantor dynamical systems naturally associated to  their renormalizations.  We now briefly discuss this approach.

An action $\Phi: \G \times \fX \to \fX$, or $(\fX, \G, \Phi)$, is said to be a  \emph{Cantor action} if  $\G$ is a finitely generated group, $\fX$ is a Cantor metric space, and the action is   minimal.    The basic properties of Cantor actions are discussed in Section~\ref{sec-prelims}.  
      In  Section~\ref{sec-renormalizable},  we associate   a minimal equicontinuous  Cantor action $(X_{\vp}, \G, \Phi_{\vp})$  to a renormalization $\vp$ of $\G$. Furthermore, we show that the  renormalization   $\vp \colon \G \to \G$ induces a \emph{renormalization   of the action} $(X_{\vp}, \G, \Phi_{\vp})$, as   in Definitions~\ref{defn-renormalizable} and \ref{defn-Vrenormalizable}.  
        
   Let $\Phi_{\vp} \colon \G \to \Homeo(X_{\vp})$ be the action homomorphism associated to the minimal equicontinuous  Cantor dynamical system  $(X_{\vp}, \G, \Phi_{\vp})$.
  The closure of the image $\Phi_{\vp}(\G)$ is an infinite 
  profinite group, denoted by  $\whGamma_{\vp}$,  called the \emph{Ellis group} of the action in the literature \cite{Auslander1988,Ellis1960,EllisGottschalk1960}. There is an induced transitive action $\whPhi_{\vp}$ of $\whGamma_{\vp}$ on $X_{\vp}$, and the \emph{discriminant group} $\cD_{\vp} \subset \whGamma_{\vp}$ is the isotropy subgroup of this action at the canonical fixed-point $x_{\vp} \in X_{\vp}$ of the contraction $\lambda_\vp \colon X_{\vp} \to X_{\vp}$. 
     The  isomorphism class of $\cD_{\vp}$ depends only on the conjugacy class of the action, and has other invariance properties   \cite{DHL2016a, DHL2016c,HL2018a,HL2019a}.  If $\G$ is abelian,   the discriminant group $\cD_{\vp}$ is   the trivial group for any renormalization, but for $\G$ non-abelian  it need not be trivial.  The regularity properties of the restricted action of $\cD_{\vp}$ on $X_{\vp}$, as discussed in   Section~\ref{sec-regularity},  play a fundamental role in the proof of  our results.
                                
  The proof of Theorem~\ref{thm-main4} is 
         given in Section~\ref{sec-finitedisc}, and uses Theorem~\ref{thm-main2} which is based on the results in Reid   \cite{Reid2014},  quoted as Theorem~\ref{thm-reid} below,   and   Theorem~\ref{thm-proaction} and Proposition~\ref{prop-embedding} in this work.  
    
      We give an example  in Section~\ref{subsec-nilpotent} of a renormalization of the Heisenberg group for which $\cD_{\vp}$ is an infinite 
  profinite group. Thus, while the assumption that $\cD_{\vp}$ is finite  is sufficient to conclude that $\G/C(\cG_{\vp})$ is nilpotent, it is not a necessary condition. On the other hand, there are renormalizations of the Heisenberg group for which $\cD_{\vp}$ is the trivial group.    
The known examples    of renormalizations  suggest  that it is   an interesting problem  to study the collection of all  renormalizations for a given group $\G$, even for the Heisenberg group.  

  We will now discuss dynamical properties of the Cantor actions  associated to renormalizations, which play the key role in the proof of Theorems \ref{thm-main4}, and of its corollary for strongly scale-invariant groups.     While most of these properties   and results do not require that $\G$ be finitely-generated,  many of our results do require this assumption, as will be pointed out when appropriate.
    
A Cantor action $(\fX, G, \Phi)$ is \emph{free} if for  any $g \in \G$ which is not the identity,   the action  $\Phi(g)$ has no fixed points. The action is \emph{topologically free}, as in  Definition~\ref{def-topfree},  if the set of points fixed by at least one non-trivial element of the group is a meager set.
The notion of a \emph{quasi-analytic} Cantor action, as in  Definition~\ref{def-qa},  was introduced in the   works \cite{DHL2016c,HL2018a}  as a generalization of the notion of a topologically free action.    
The quasi-analytic property of a Cantor action is a fundamental property of  renormalizable groups and actions. 

 \begin{thm}\label{thm-main1}
 The  Cantor action  $\Phi_{\vp} \colon \G \times X_{\vp} \to X_{\vp}$ associated to a renormalization $\vp$ is quasi-analytic.  Hence, if the action $\Phi_{\vp}$ is   effective  then it is topologically free.   
  \end{thm}
 
  In fact, Theorem \ref{thm-main1} is a consequence of a stronger statement.
Given a Cantor action  $(\fX,\G,\Phi)$, let  $\Phi(\G) \subset \Homeo(\fX)$ denote the image subgroup. If the action is equicontinuous, then  the  closure   $\overline{\Phi(\G)} \subset \Homeo(\fX)$ in the \emph{uniform topology of maps} is a separable profinite group. This is discussed further   in Section~\ref{subsec-models}. 
For the Cantor action 
$(X_{\vp}, \G, \Phi_{\vp})$ associated to a renormalization $\vp$, we denote this closure by  $\whGamma_{\vp} = \overline{\Phi_{\vp}(\G)}$.  We prove in  Theorem~\ref{thm-proaction}   that the profinite action $\whPhi_{\vp} \colon \whGamma_{\vp} \times X_{\vp} \to X_{\vp}$ is quasi-analytic, which implies Theorem~\ref{thm-main1}. The quasi-analytic property is used   to prove that the monomorphism $\vp: \G \to \G$ induces an open embedding of the closure $\whGamma_{\vp}$ into itself as below:   
  
 \begin{thm}\label{thm-embedding}
    Let $\vp$ be a renormalization of the finitely generated group $\G$. Then $\vp$ induces an injective  homomorphism  $\whvp \colon \whGamma_{\vp} \to \whGamma_{\vp}$ whose image is a clopen subgroup of $\whGamma_{\vp}$. 
  \end{thm}
This is proved in Section~\ref{sec-open}, where we use its conclusion to  obtain a structure theorem for the closure $\whGamma_\varphi$, an important consequence of which is the following: 
 \begin{thm}\label{thm-dynamicdiscriminant}
    Let $\vp$ be a renormalization of the finitely generated group $\G$, and     $\whvp \colon \whGamma_{\vp} \to \whGamma_{\vp}$ the induced map given by Theorem~\ref{thm-embedding}. Then
 $\ds \cD_{\vp} =   \bigcap_{n > 0} \ \whvp^n(\whGamma_{\vp})$.
  \end{thm}
  Theorems~\ref{thm-main4} and \ref{thm-dynamicdiscriminant} yield an answer to the profinite version of the Nekrashevych-Pete question:

  \begin{thm}\label{thm-contractionprinciple}
  Let $\vp$ be a renormalization of the finitely generated group $\G$. Suppose that  
$$
K(\cG_{\vp}) = \bigcap_{\ell > 0} \vp^{\ell}(\G)  \subset \G \quad , \quad  \cD_{\vp} =   \bigcap_{n > 0} \ \whvp_0^n(\whGamma_{\vp}) \subset \whGamma_{\vp}
$$
are both finite groups.  Then $\G$ is  virtually nilpotent. If both groups  are trivial, then $\G$ is nilpotent.
  \end{thm}

 As mentioned above, our study of renormalizable groups  naturally suggests a related notion, that of   a \emph{renormalizable} equicontinuous Cantor action, as   introduced  in Definition~\ref{defn-renormalizable}. It  is modeled on the concept of a renormalizable dynamical system, and also that of  self-similar groups \cite{Grigorchuk2011,Nekrashevych2005} and percolation theory \cite[Section~9.2]{Sapir2007}. 
  The class of renormalizable Cantor actions includes the class of actions associated to a self-embedding $\varphi: \G \to \G$, discussed above. We   introduce a variant of this notion, the \emph{virtually renormalizable} actions, in Definition~\ref{defn-Vrenormalizable}. 
 
The study of renormalizable  Cantor actions is motivated, in part,  by the examples of Cantor actions on $d$-adic trees, where elements of the acting group are defined recursively, in terms of the  action of a finite set of generators on a rooted $d$-adic tree for  $d \geq 2$, where there is an embedding   $\vp \colon \G \to \G$ whose image is a subgroup of the stabilizer group of a branch of a tree (see for example \cite{Nekrashevych2005}). 
The image $\vp(\G) \subset \G$ need not be of finite index in $\G$, even though the stabilizer group of a branch always has finite index in $\G$.   We discuss the relation of the notion of a renormalizable action with the notion of self-replicating groups in more detail  in Section~\ref{subsec-renormalizable}.

  An equicontinuous Cantor action  can be \emph{quasi-analytic}, \emph{locally quasi-analytic} or \emph{wild}, using the notions and classification introduced  in the works \cite{HL2018a,HL2018b,HL2019a}.  For renormalizable actions, there is the following dichotomy:
\begin{thm}\label{thm-dichotomy}
A renormalizable equicontinuous Cantor action $(\fX,\G,\Phi)$ is either quasi-analytic, 
and in this case  $\G$ is renormalizable  and the action is topologically conjugate to the action given by a renormalization of $\G$, or the action is wild.
 \end{thm}
 This result motivates the study of  the invariants of renormalizable Cantor actions, both to understand the  invariants of the renormalization map, and to discover invariants of these actions which   distinguish between the  quasi-analytic and wild cases of Theorem~\ref{thm-dichotomy}.  Our final results in this work considers their invariant properties under continuous orbit equivalence.

\begin{thm}\label{thm-renormal}
Let $(\fX, \G, \Phi)$ and $(\fX', \G', \Phi')$ be minimal equicontinuous Cantor actions which are continuously orbit equivalent. 
If   $(\fX, \G, \Phi)$ is   renormalizable and locally quasi-analytic, then  $(\fX', \G', \Phi')$ is  virtually renormalizable.
\end{thm}

As a consequence,  we obtain that the isomorphism class of the discriminant group $\cD_{\vp}$ associated to a renormalization $\vp$   is an invariant of continuous orbit equivalence.
   \begin{thm}\label{thm-Dcoe}
Let $(X_{\vp}, \G, \Phi_{\vp})$ and   $(X'_{\vp'}, \G', \Phi'_{\vp'})$ be     Cantor actions associated to renormalizations 
$\vp \colon \G \to \G$ and $\vp' \colon \G' \to \G'$, respectively. If the actions are   continuously orbit equivalent, 
  then  their discriminant groups $\cD_{\vp}$  and $\cD'_{\vp'}$ are isomorphic. 
  \end{thm}

 Examples and applications of our results   are discussed in Section~\ref{sec-examples}.
       
 Section~\ref{sec-problems} discusses     open problems. In particular, 
the works \cite{ALBLLN2020,HL2018a,HL2019a} study wild Cantor actions, and the relations between the discriminant invariant and the wild property for the action.  It is an interesting problem to further explore this relation for renormalizable actions, as these include many class of branch groups and related constructions, as in \cite{Bartholdi2003,BGS2012,Nekrashevych2005,NekkyPete2011,Nekrashevych2018,GL2019}. 

\subsection*{Acknowledgments} We would like to thank an anonymous referee for helpful comments and suggestions. OL is supported by FWF Project P31950-N35. During this work, WvL was supported by NSF Grant DMS-1855371, the Max Planck Institute for Mathematics, and NSF Grant DMS-1928930 while participating in a program hosted
by the Mathematical Sciences Research Institute in Berkeley, California during the Fall 2020 semester.

\section{Cantor actions}\label{sec-prelims}

In this section, we recall some of the 
 properties of     Cantor actions.
 A basic reference  is   \cite{Auslander1988}.

\subsection{Basic concepts}\label{subsec-basics}
For an action    $\Phi \colon \G \times \fX \to \fX$    and $\gamma \in \G$,  let  $\gamma   x = \Phi(\gamma)(x)$. We also sometimes write $\gamma \cdot x = \Phi(\gamma)(x)$ when necessary for notational clarity.

Let  $(\fX,\G,\Phi)$   denote an action  $\Phi \colon \G \times \fX \to \fX$.
The orbit of  $x \in \fX$ is the subset $\cO(x) = \{\gamma   x \mid \gamma \in \G\}$. 
The action is \emph{minimal} if  for all $x \in \fX$, its   orbit $\cO(x)$ is dense in $\fX$.

An action  $(\fX,\G,\Phi)$ is \emph{equicontinuous} with respect to a metric $\dX$ on $\fX$, if for all $\e >0$ there exists $\delta > 0$, such that for all $x , y \in \fX$   with 
 $\ds  \dX(x,y) < \delta$ and  all  $\gamma \in \G$,    we have  $\dX(\gamma   x, \gamma   y) < \e$.
The equicontinuous property   is independent of the choice of the metric   on $\fX$  by Proposition \ref{prop-CO}.

 An action $(\fX,\G,\Phi)$  is \emph{effective}, or \emph{faithful},  if the  action homomorphism $\Phi \colon \G \to \Homeo(\fX)$ has  trivial kernel.
 The  action  is \emph{free} if for all $x \in \fX$ and $\gamma \in \G$,   $\gamma x = x$ implies that $\gamma = e$,   the identity of the group. The \emph{isotropy group} of $x \in \fX$ is the subgroup 
\begin{equation}\label{eq-isotropyx}
\G_x = \{ \gamma \in \G \mid \gamma x = x\} \ . 
\end{equation}

Let $\Fix(\gamma) = \{x \in \fX \mid \gamma x = x \}$, and introduce the \emph{isotropy set}
\begin{equation}\label{eq-isotropy}
 \Iso(\Phi) = \{ x \in \fX \mid \exists ~ \gamma \in \G ~ , ~ \gamma \ne id ~, ~\gamma x = x    \} = \bigcup_{e \ne \gamma \in \G} \ \Fix(\gamma) \ . 
\end{equation}

  \begin{defn}\cite{BoyleTomiyama1998,Li2018,Renault2008} \label{def-topfree}
  $(\fX,\G,\Phi)$ is said to be \emph{topologically free}  if  $\Iso(\Phi) $ is meager in $\fX$. 
 \end{defn}
 
Note that if $\Iso(\Phi)$ is meager, then $\Iso(\Phi)$ has empty interior. 
 
   The   notion of topologically free Cantor actions was introduced by Boyle in his thesis \cite{Boyle1983}, and later  used  in the works by Boyle and Tomiyama \cite{BoyleTomiyama1998}   for the study of classification of   Cantor actions,   by   Renault \cite{Renault2008}
     for the study of the $C^*$-algebras associated to Cantor actions, and by   Li \cite{Li2018} in his study of    rigidity properties of equicontinuous Cantor actions.

Now assume that $\fX$ is a Cantor space. 
Let $\CO(\fX)$ denote the collection  of all clopen (closed and open) subsets of  $\fX$, which forms a basis for the topology of $\fX$. 
For $\phi \in \Homeo(\fX)$ and    $U \in \CO(\fX)$, the image $\phi(U) \in \CO(\fX)$.  
The following   result is folklore, and a proof is given in \cite[Proposition~3.1]{HL2018b}.
 \begin{prop}\label{prop-CO}
 A Cantor  action   $\Phi \colon \G \times \fX \to \fX$  is  equicontinuous  if and only if  the orbit of every $U \in \CO(\fX)$ is finite for the induced action $\Phi_* \colon \G \times \CO(\fX) \to \CO(\fX)$.
\end{prop}
 
 Let  $(\fX,\G,\Phi)$ be a minimal equicontinuous Cantor action.  
We say that $U \subset \fX$  is \emph{adapted} to the action   if $U$ is a   non-empty clopen subset, and for any $\gamma \in \G$, 
if $\Phi(\gamma)(U) \cap U \ne \emptyset$ then  $\Phi(\gamma)(U) = U$.   The proof of  Proposition~3.1 in \cite{HL2018b} shows that given  $x \in \fX$ and a   clopen set $W$ with $x \in W$, there is an adapted clopen set $U$ with $x \in U \subset W$. 

The key property of an   adapted set $U$  is that the set of ``return times'' to $U$, 
 \begin{equation}\label{eq-adapted}
\G_U = \left\{\gamma \in \G \mid \Phi(\gamma)(U) = U   \right\}  
\end{equation}
is a subgroup of   $\G$, called the \emph{stabilizer} of $U$.      
  Then for $\gamma, \gamma' \in \G$ with $\Phi(\gamma)(U) \cap \Phi(\gamma')(U) \ne \emptyset$ we have $\Phi(\gamma^{-1}) \circ \Phi(\gamma')(U) = U$, hence $\gamma^{-1} \ \gamma' \in \G_U$. Thus, as the action is assumed to be minimal,  the  translates $\{ \Phi(\gamma)(U)  \mid \gamma \in \G\}$ form a finite clopen partition of $\fX$, and are in 1-1 correspondence with  the elements in the quotient space $X_U = \G/\G_U$. Then $\G$ acts by permutations of the finite set $X_U$ and so the stabilizer group $\G_U \subset \G$ has finite index.  
  
The action of $\gamma \in \G$ on $X_U$ is trivial precisely when $\gamma$ is a stabilizer of each coset $\kappa \cdot \G_U$, so   $\gamma \in C_U$ where
$C_U \ = \  \bigcap_{\kappa \in \G} \ \kappa  \ \G_U \ \kappa^{-1} \ \subset \ \G_U$ 
  is the largest normal subgroup of $\G$ contained in $\G_U$. 
The action of    the finite group $Q_U \equiv \G/C_U$ on $X_U$ by permutations is a finite approximation of the action of $\G$ on $\fX$,  and the isotropy group of the identity coset   $e \cdot \G_U$ is  $D_U \equiv \G_U/C_U \subset Q_U$. 

\begin{defn}\label{def-adaptednbhds}
Let  $(\fX,\G,\Phi)$   be an equicontinuous  Cantor    action.
A properly descending chain of clopen sets $\cU = \{U_{\ell} \subset \fX  \mid \ell \geq 0\}$ is   an \emph{adapted neighborhood basis} at $x \in \fX$ for the action $\Phi$,   if
    $x \in U_{\ell +1} \subset U_{\ell}$ for all $ \ell \geq 0$,        each $U_{\ell}$ is adapted to the action $\Phi$, and the intersection      $\cap_{\ell > 0}  \ U_{\ell} = \{x\}$.
\end{defn}
Given $x \in \fX$ and   $\e > 0$, Proposition~\ref{prop-CO} implies there exists an adapted clopen set $U \in \CO(\fX)$ with $x \in U$ and $\diam(U) < \e$.  Thus, one can choose a descending chain $\cU$ of adapted sets in $\CO(\fX)$ whose intersection is $x$, which shows the following:

\begin{prop}\label{prop-adpatedchain}
Let  $(\fX,\G,\Phi)$   be a minimal equicontinuous  Cantor    action. Given $x \in \fX$, there exists an adapted neighborhood basis $\cU$ at $x$ for the action $\Phi$.
 \end{prop}

\subsection{The dynamical profinite model}\label{subsec-models}

Given an equicontinuous Cantor action  $(\fX,\G,\Phi)$, let  $\Phi(\G) \subset \Homeo(\fX)$ denote the image subgroup. Then the  closure $\overline{\Phi(\G)} \subset \Homeo(\fX)$ in the \emph{uniform topology of maps} is a separable profinite group. This group is identified with  the \emph{Ellis group} for the action, as defined in \cite{Auslander1988,Ellis1960,EllisGottschalk1960}; see also \cite[Section~2]{DHL2016a}. 
Each  element   $\whg \in \overline{\Phi(\G)}$ is the uniform limit of a sequence of maps   $\{\Phi(g_i) \mid i \geq 1\} \subset \Phi(\G)$. We sometimes       denote   the limit $\whg$ by $(g_i)$. 
 
For example, if $\G$ is an abelian group, then $\overline{\Phi(\G)}$ is a compact totally disconnected abelian group, which can be thought of as the group of asymptotic motions of the system. When $\G$ is non-abelian,    the   closure $\overline{\Phi(\G)}$ can have  much more subtle algebraic properties.  
 
Let $\whPhi \colon \overline{\Phi(\G)} \times \fX \to \fX$ denote the induced action of $\overline{\Phi(\G)}$ on $\fX$. For $\whg \in \overline{\Phi(\G)}$ we   write its action on $\fX$ by $\whg \, x = \whPhi(\whg)(x)$.
For a minimal equicontinuous  action $\Phi \colon \G \times \fX \to \fX$, the group $\overline{\Phi(\G)}$ acts transitively on $\fX$.    
Given $x \in \fX$,   introduce the isotropy group  at $x$,  
\begin{align}\label{iso-defn2}
 \overline{\Phi(\G)}_x = \{ \whg  \in \overline{\Phi(\G)} \mid \whg \, x = x\} \subset \Homeo(\fX) \ ,
\end{align}
which is a closed subgroup of $\overline{\Phi(\G)}$, and    thus is either finite, or is an infinite profinite group.    

\begin{defn}\label{def-discriminant}
The group $\overline{\Phi(\G)}_x$ is called the \emph{discriminant} of the action $(\fX,\G,\Phi)$.
\end{defn}

There is a   natural identification $\fX \cong \overline{\Phi(\G)}/\overline{\Phi(\G)}_x$ of left $\overline{\Phi(\G)}$-spaces, and thus      
the conjugacy class of $\overline{\Phi(\G)}_x$ in $\overline{\Phi(\G)}$ is independent of the choice of   $x$.    If $\overline{\Phi(\G)}_x$ is the trivial group, then $\fX$ is identified with the profinite group $\overline{\Phi(\G)}$, and the   action $\whPhi$ is  free.  Note that  there exists examples of   free minimal equicontinuous Cantor actions $(\fX,\G,\Phi)$ for which the action $\whPhi$ is not free, and in fact  $\overline{\Phi(\G)}_x$ is an  infinite profinite group. The first such examples were constructed by  Fokkink and Oversteegen   in \cite[Section~8]{FO2002}, and further examples are constructed in   \cite[Section~10]{DHL2016c}.
  
\subsection{Equivalence of Cantor actions}\label{subsec-equivalence}

We recall three notions of equivalence of  Cantor actions.

 The first and strongest   is that  of 
  {isomorphism} of Cantor actions, which   is a   generalization  of the usual notion of conjugacy of topological actions. For $\G = \mZ$, isomorphism corresponds to the notion of ``flip conjugacy'' introduced    in the work of Boyle and Tomiyama \cite{BoyleTomiyama1998}.
  
 \begin{defn} \label{def-isomorphism}
Two Cantor actions $(\fX, \G, \Phi)$ and $(\fX', \G', \Phi')$   are said to be \emph{isomorphic}  if there is a homeomorphism $h \colon \fX \to \fX'$ and a group isomorphism $\Theta \colon \G \to \G'$ so that 
\begin{equation}\label{eq-isomorphism}
\Phi(\gamma) = h^{-1} \circ \Phi'(\Theta(\gamma)) \circ h   \in   \Homeo(\fX') \  {\rm for \ all} \ \gamma \in \G \ .
\end{equation}
 \end{defn}

 \emph{Return equivalence}  is a form of ``virtual   isomorphism'' for minimal equicontinuous Cantor actions, and  is   weaker than the notion of isomorphism. This equivalence  is natural when considering the Cantor systems arising from geometric constructions, as    in   the works  \cite{HL2018a,HL2018b,HL2019a}.

 Throughout this work, by a small abuse of  notation, for a minimal equicontinuous Cantor action  $(\fX,\G,\Phi)$ and adapted subset $U \subset \fX$, we use $\Phi_U$ to denote both the restricted action $\Phi_U \colon \G_U \times U \to U$ and the induced quotient action $\Phi_U \colon H_U \times U \to U$ where $H_U = \Phi(\G_U) \subset \Homeo(U)$.

   \begin{defn}\label{def-return}
Minimal equicontinuous Cantor  actions $(\fX, \G, \Phi)$ and $(\fX', \G', \Phi')$ are  \emph{return equivalent} if there exists 
  an adapted set $U \subset \fX$ for the action $\Phi$   and  
  an adapted set $V \subset \fX'$ for the action $\Phi'$,
such that   the  restricted actions $(U, H_U, \Phi_U)$ and $(V, H'_V, \Phi'_V)$ are isomorphic.
\end{defn}

  \emph{Continuous orbit equivalence} for   Cantor actions  was introduced in     \cite{Boyle1983,BoyleTomiyama1998}, and plays a fundamental role in various approaches to the classification of these actions \cite{Renault2008}.
Consider the equivalence relation on $\fX$ defined by an action    $(\fX,\G,\Phi)$,
 \begin{equation}\label{eq-ERX}
\cR(\fX, \G, \Phi) \equiv \{(x,  \gamma x)) \mid x \in \fX, \gamma \in \G\} \subset \fX \times \fX ~.
\end{equation}
Given   actions $(\fX,\G,\Phi)$ and $(\fX', \G', \Phi')$, we say they are \emph{orbit equivalent} if there exist a bijection   $h \colon \fX \to \fX'$ which maps 
$\cR(\fX, \G, \Phi)$   onto   $\cR(\fX', \G', \Phi')$, and similarly for the inverse map $h^{-1}$.

\begin{defn}\label{def-torb1}
Let $(\fX,\G,\Phi)$ and $(\fX', \G', \Phi')$ be      Cantor actions.
A \emph{continuous orbit equivalence}  between the actions  is a homeomorphism    $h \colon \fX \to \fX'$ which is an orbit equivalence, and there exist continuous functions $\alpha \colon \G \times \fX \to \G'$ and $\beta \colon \G' \times \fX' \to \G$  such that:
\begin{enumerate}
\item  for each $x \in \fX$ and $\gamma \in \G$, there exists $\alpha(\gamma,x) \in \G'$ and an open set $x \in U_x \subset \fX$ such that 
$ \Psi(\alpha(\gamma,x))  \circ h | U_x = h \circ \Phi(\gamma) | U_x$; \\
\item for each $y \in \fX'$ and $\gamma' \in \G'$, there exists $\beta(\gamma',y) \in \G$ and an open set $y \in V_y \subset \fX'$ such that 
$ \Phi(\beta(\gamma',y))  \circ h | V_y = h \circ \Psi(\gamma') | V_y$.
\end{enumerate}
\end{defn}
Note that   the maps $\alpha$ and $\beta$ are not assumed to be cocycles over the respective actions.

\section{Renormalizable groups}\label{sec-renormalizable}

 In this section, we construct the minimal equicontinuous Cantor action $(X_{\vp}, \G, \Phi_{\vp})$ associated to a renormalization $\vp \colon \G \to \G$, and give some of the basic properties of this action.

Set $\G_0 = \G$, and   for $\ell \geq 1$,  recursively define  subgroups  
 $\G_{\ell} \subset \G$, where $\G_{\ell}  = \vp(\G_{\ell-1}) \equiv \vp^{\ell}(\G)$.

   Let $\cG_{\vp} \equiv \{\G_{\ell} \mid \ell \geq 0\}$ denote the   descending group chain, where each $\G_{\ell}$ has finite index in $\G$. 
Denote the intersection of the group chain by $K(\cG_{\vp}) \equiv \bigcap_{\ell > 0} \ \G_{\ell}$.
If $K(\cG_{\vp})$  is a finite group, then the group $\G$  is said to be \emph{strongly scale-invariant},   in the terminology of Nekrashevych and Pete \cite{NekkyPete2011}.

Let $X_{\ell} = \G/\G_{\ell}$ be the finite    coset space. Note that $X_{\ell}$   is not necessarily a group, as the subgroup $\G_{\ell}$ is not assumed to be normal in $\G$. 
   Note that  $\G$ acts transitively on the left on   $X_{\ell}$, and  the inclusion $\G_{\ell +1} \subset \G_{\ell}$ induces a natural $\G$-invariant quotient map $p_{\ell +1} \colon X_{\ell +1} \to X_{\ell}$.
The inverse limit space
 \begin{equation} \label{eq-invlimspace}
X_{\vp} \equiv \varprojlim \ \{p_{\ell +1} \colon X_{\ell +1} \to X_{\ell} \mid \ell \geq 0\}   
\end{equation}
with the Tychonoff topology  is a Cantor space.
The actions of $\G$ on the factors $X_{\ell}$ induce    a minimal  equicontinuous action on $X_{\vp}$, 
denoted by  $\Phi_{\vp} \colon \G \times X_{\vp} \to X_{\vp}$ or by $(X_{\vp}, \G, \Phi_{\vp})$.  

  Let $\whGamma_{\vp} \subset \Homeo(X_{\vp})$ denote the closure of the image $\Phi_{\vp}(\G) \subset \Homeo(X_{\vp})$,  as introduced in Section~\ref{subsec-models}, and let $\whPhi_{\vp} \colon \whGamma_{\vp} \times X_{\vp} \to X_{\vp}$ denote the induced action.

The embedding $\vp$ induces a  mapping   $\lambda_{\vp} \colon X_{\vp}  \to X_{\vp}$, which is defined as the shift map on   sequences as follows. A point $\whx \in X_{\vp}$ is defined by an equivalence class of sequences $\whx = (g_0, g_1, g_2, \ldots)$ with each $g_{\ell} \in \G$ satisfying the relations  
$g_{\ell}  = g_{\ell +1} \ {\rm mod} \  \G_{\ell}$ for all $\ell \geq 0$. Then 
$\lambda_{\vp}(\whx) =   (e, \vp(g_0), \vp(g_1), \vp(g_2), \ldots)$ is well-defined, and is a contraction  on $X_{\vp}$.
Let $x_{\vp} \in X_{\vp}$ be the unique fixed point for  $\lambda_{\vp}$. 
Then  $x_{\vp} = (e,e,e,\ldots )$ where $e \in \G$ is the identity.

 Following Definition~\ref{def-discriminant}, we obtain   a fundamental notion associated to a renormalization of $\G$. 
    \begin{defn}\label{def-vpdiscriminant}
   The \emph{discriminant group}   of    $(X_{\vp}, \G, \Phi_{\vp})$
   is   $\cD_{\vp} \equiv   (\whGamma_{\vp})_{x_{\vp}} \subset \whGamma_{\vp}$.
     \end{defn}

 For  $k \geq 0$, define 
  \begin{equation} \label{eq-clopenbasis}
  U_k = \{(g_0,g_1,g_2,\ldots) \in X_{\vp} \mid g_i=e \textrm{ for } 0 \leq i \leq k \} \cong \varprojlim \ \{p_{\ell +1} \colon \G_k/\G_{\ell +1} \to \G_k/\G_{\ell} \mid \ell \geq k\} ~ , 
\end{equation}
which is a clopen subset of $X_{\vp}$ adapted to the action $\Phi_{\vp}$, with    stabilizer subgroup     $\G_{U_k} = \G_k$.
The  clopen sets  $\{U_k \mid k \geq 0\}$ form an adapted neighborhood basis at $x_{\vp}$, and so we have $x_{\vp}   = \bigcap_{k > 0} \ U_k$.  

Observe that   for all $\ell \geq 0$, the contraction mapping $\lambda_{\vp} \colon X_{\vp} \to X_{\vp}$ defined above restricts to a homeomorphism  onto $\lambda_{\vp} \colon U_{\ell} \to U_{\ell +1}$.
  
 As the orbit of $x_{\vp}$ is dense in $X_{\vp}$,  for any non-empty open subset $U \subset X$ there exists $g \in \G$ so that $\Phi(g)(x_{\vp}) \in U$. It follows that there also exists  $k > 0$ such that $\Phi(g)(U_k) \subset U$.
 
\subsection{The geometric (tree) model}\label{subsec-treemodels}

Given a group chain $\cG_{\vp}$ as above, the associated dynamical system $\Phi_{\vp} \colon \G \times X_{\vp} \to X_{\vp}$ can be represented as an action of a subgroup of the automorphism group of a regular rooted tree as we discuss now. The construction we discuss is applicable to any chain $\G_0 \supset \G_1 \supset \G_2 \supset$ of finite index subgroups of the group $\G=\G_0$, without a requirement that the subgroups in the chain are isomorphic to $\G$.

A \emph{tree} $T$ consists of a set of vertices $V = \bigsqcup_{\ell \geq 0} V_\ell$, where $V_\ell$ is a finite vertex set at level $\ell$, and of edges joining vertices in $V_{\ell+1}$ and $V_\ell$, for all $\ell \geq 0$, defined as follows. For $\ell \geq 0$, identify the vertex set $V_\ell$ with the coset space $X_\ell = \G/\G_\ell$. Join $v_\ell \in V_\ell$ and $v_{\ell+1} \in V_{\ell+1}$ by an edge if and only if $v_{\ell+1} \subset v_\ell$ as cosets. Let $d = |\G : \vp(\G)|$ be the index, then for $\ell \geq 0$ the cardinality of $V_\ell$ is $d^\ell$. Such a tree is called a $d$-ary, or a regular tree.

An infinite path in $T$ is a sequence of vertices $(v_\ell)_{\ell \geq 0}$ such that $v_{\ell+1}$ and $v_\ell$ are joined by an edge, for $\ell \geq 0$. The boundary $\partial T$ of $T$ is the collection of all infinite paths in $T$, and so it is the subspace
  $$\partial T = \{(v_\ell)_{\ell \geq 0} \subset \prod_{\ell \geq 0} V_\ell \mid v_{\ell+1} \textrm{ and }v_\ell \textrm{ are joined by an edge}\}.$$
The space $\partial T$ is a Cantor set with the relative topology from the product topology on $ \prod_{\ell \geq 0} V_\ell $. It is immediate that the identification of the vertex sets $V_\ell$ with the coset spaces $X_\ell$ induces an identification of $\partial T$ with the inverse limit space $X_\vp$ defined by \eqref{eq-invlimspace}, with points in $X_\vp$ corresponding to infinite paths in $\partial T$. 

The action of $\G$ on the coset spaces $X_\ell=V_\ell$, $\ell \geq 0$, is by permutations. Since the action of $\G$ preserves the containment of cosets, the action preserves the connectedness of the tree $T$, that is, the vertices $v_{\ell} \in V_{\ell}$ and $v_{\ell+1} \in V_{\ell+1}$ are joined by an edge if and only if for any $g \in G$ the images $g \cdot v_\ell \in V_{\ell}$ and $g \cdot v_{\ell+1} \in V_{\ell+1}$ are joined by an edge. Thus every $g \in \G$ defines an automorphism of the tree $T$, and we can consider $\G$ as a subgroup of the group of tree automorphisms $Aut(T)$. 

The study of actions of subgroups of automorphism groups of regular trees is an active topic in geometric group theory, see for instance \cite{Nekrashevych2005}. When studying the dynamical properties of an action $(X_\vp, \G, \Phi_{\vp})$, sometimes it is useful to represent it as an action on the boundary of a regular tree. However, our results in this paper rely heavily on the properties of the profinite completion $\fG$ of $\G$, and the combinatorial character of the methods used to study group actions on rooted trees makes their use in the study of profinite completions rather cumbersome. In this paper, we rely mostly on the algebraic methods we describe further in Sections \ref{subsec-algmodels}  and \ref{subsec-universalmodel}, while explaining the implications of our results for the actions of automorphisms of regular rooted trees in Section~\ref{subsec-renormalizable}.

\subsection{The algebraic profinite model}\label{subsec-algmodels} 
We next introduce an alternate profinite model for the minimal equicontinuous  action $(X_{\vp}, \G, \Phi_{\vp})$.
For each $\ell \geq 1$, let $C_{\ell}$ denote the largest normal subgroup (the \emph{core}) of the stabilizer group $\G_{\ell}$, so
\begin{equation}\label{eq-core}
C_{\ell} ~ = ~ \bigcap_{g \in \G} ~ g \ \G_{\ell} \ g^{-1} ~ \subset \G_{\ell} ~ .
\end{equation}
As   $\G_{\ell}$ has finite index in $\G$, the same holds for $C_{\ell}$. Observe that for all $\ell \geq 1$,   we have $C_{\ell +1} \subset C_{\ell}$.
Introduce the quotient group  $Q_{\ell} = \G/C_{\ell}$ with identity element $e_{\ell} \in Q_{\ell}$. There are natural quotient maps $q_{\ell+1} \colon Q_{\ell +1} \to Q_{\ell}$, and we can form the inverse limit group
\begin{equation}\label{eq-invgroup}
\whGamma_{\infty} = \varprojlim \ \{q_{\ell +1} \colon Q_{\ell +1} \to Q_{\ell} \mid \ell \geq 0\} ~ .
\end{equation}
\begin{thm}\cite[Theorem~4.4]{DHL2016a}
There is a natural isomorphism $\whtau \colon \whGamma_{\vp} \to \whGamma_{\infty}$ which identifies the discriminant group $\cD_{\vp}$ with the inverse limit group
\begin{equation}\label{eq-discformula}
\cD_{\infty}  = \varprojlim \ \{q_{\ell +1} \colon \G_{\ell +1}/C_{\ell+1} \to \G_{\ell}/C_{\ell} \mid \ell \geq 0\} \subset \whGamma_{\infty}~ .
\end{equation}
\end{thm}

 There is an interpretation of the group $\cD_{\infty}$  as an \emph{asymptotic defect} of the $\G$-action on $X_{\infty}$.  Suppose that $\G_{\ell}$ is a normal subgroup, so that the quotient  $\G/\G_{\ell}$ is a group. Then $\G/\G_{\ell}$ acts transitively and freely on $X_{\ell}$.  For example, if $\G$ is abelian then this is always true. In general, for the normal core $C_{\ell} \subset \G_{\ell}$, the finite group $Q_{\ell} = \G/C_{\ell}$ acts transitively on $X_{\ell}$ and the finite subgroup $D_{\ell} = \G_{\ell}/C_{\ell}$ is the ``defect''    for the action of $Q_{\ell}$ on $X_{\ell}$ being a free action.  Then $\cD_{\infty}$ is the inverse limit of these finite defects, and provides a measure of the deviation of the action  $\whPhi_{\infty}$ of $\whGamma_{\infty}$ on $X_{\infty}$ from being a free action.

 Associated to the group chain $\cG_{\vp}$,  there are two subgroups,
\begin{equation}\label{eq-kernels}
K(\cG_{\vp}) = \bigcap_{\ell > 0} \G_{\ell} \quad , \quad C(\cG_{\vp}) = \bigcap_{g \in \G} g \ K(\cG_{\vp}) \ g^{-1}~ .
\end{equation}
where     $\ds C(\cG_{\vp})$ is the largest normal subgroup of  $\G$ contained in $\ds K(\cG_{\vp})$.
Note that  for any $g \in C(\cG_{\vp})$, the action of $\Phi_{\vp}(g)$ on $X_{\vp}$ is trivial.

\subsection{The universal profinite model}\label{subsec-universalmodel}

We now introduce a  model for the action  $(X_{\vp}, \G, \Phi_{\vp})$ in terms of the profinite completion $\fG$ of $\G$.
 Recall that $\fG$ is the inverse limit of the finite quotient groups $\G/N$,  for the   set 
$\ds \cN = \{N    \mid N \subset \G ~ \mbox{is a normal subgroup of finite index}\}$ which is partially ordered by inclusion. That is,   
$\ds \fG = \varprojlim ~ \{\G/N \mid N \in \cN\}$.

There is a homomorphism $\psi \colon \G \to \fG$ with dense image, and the kernel of $\psi$ is the group $N({\psi})$ given by the intersection of all normal subgroups of finite index in $\G$.  Thus, $N({\psi})$ is trivial exactly when the group $\G$ is residually finite. 

By the universal property of the profinite completion, the    map $\Phi_{\vp} \colon \G \to \whGamma_{\vp} \subset \Homeo(X_{\vp})$ induces a surjective map $\Pi_{\vp} \colon \fG \to \whGamma_{\vp}$ of profinite groups,  
  and let    $\fN(\Pi_{\vp}) \subset \fG$ denote its  kernel.

 Let $\fD_{\vp} \equiv \fG_{x_{\vp}} \subset \fG$ denote the isotropy subgroup   at $x_{\vp}$ of the action $\whPhi_{\fG}$ of  $\fG$.
Then $\fN(\Pi_{\vp}) \subset \fD_{\vp}$.  

 We use the universal property of $\fG$ to show a  basic fact   required for our study of renormalizations.
\begin{prop}\label{prop-profiniteG}
The renormalization $\vp \colon \G \to \G$ induces an    open embedding $\whvp_{\fG} \colon \fG \to \fG$.
\end{prop}
\proof
Let $\G_1 =\vp(\G) \subset \G$ denote the image of $\vp$. Then the partially ordered set
$$\cN_1 = \{N_1    \mid N_1 \subset \G_1 ~ \mbox{is a   subgroup of finite index, normal in $\G_1$}\}$$
yields the universal profinite completion $\fG_1$ of $\G_1$ and 
  $\vp$ induces an   isomorphism $\whvp_1 \colon \fG \to \fG_1$.

Next, consider the partially ordered set 
$\cN_1' = \{N \cap \G_1    \mid N \in \cN\}$, where $\cN$ is the collection of finite-index   subgroups which are normal in  $\G$. 
Let $\fG_1'$ denote the profinite completion of $\G_1$ with respect to $\cN_1'$. It is immediate that   $\fG_1'$ is the closure of $\psi(\G_1) \subset \fG$ in $\fG$.

 Note that if $N' \in \cN'$, then $N'$ is also normal in $\G_1$ so $N' \in \cN_1$. Conversely, if $N \in \cN_1$ then its normal core $CN_1 = \cap_{g \in \G} \ g^{-1} N_1 \ g \subset N_1$ is a normal subgroup of $\G$ which has finite index in $\G_1$ so $CN_1 \in \cN_1'$. 
Thus, the two partially ordered sets $\cN_1$ and $\cN_1'$ are cofinal in $\G_1$, hence the identity map induces a homeomorphism  $\whId \colon \fG_1 \cong \fG_1'$. It follows that the composition $\whvp_{\fG} \equiv \whId \circ \whvp_1 \colon \fG \to \fG_1 \cong \fG_1'$ is an isomorphism onto the subgroup   $\fG_1' \subset \fG$ of finite index, which is thus open.
\endproof

For each $\ell \geq 1$, let $\fG_{\ell} = \whvp_{\fG}^{\ell}(\fG)$ which is   an open   subgroup  of $\fG$ of finite index, with $\fG/\fG_{\ell} \cong \G/\G_{\ell}$. 

By an argument analogous to the proof of Proposition~\ref{prop-profiniteG}, $\fG_{\ell}$ is identified with the closure of $\G_{\ell}$ in $\fG$, and thus $\fG_{\ell} = \{ g \in \fG \mid \whPhi_{\fG}(g)(U_{\ell}) = U_{\ell} \}$. That is, $U_{\ell}$ is an adapted set for 
the action $\whPhi_{\fG} \colon \fG \times X_{\vp} \to X_{\vp}$ with stabilizer $\fG_{U_{\ell}} = \fG_{\ell}$.
It follows that the isotropy group $\fD_{\vp}$ at $x_{\vp} \in X_{\vp}$ of the action $\whvp_{\fG}$   is given by   $  \cap_{\ell\geq 0} \ \fG_{\ell}$, and so $\Pi_{\vp}(\fD_{\vp}) = \cD_{\vp} \subset \whGamma_{\vp}$.

\begin{cor}\label{cor-automorphism}
The embedding $\whvp_{\fG} \colon \fG \to \fG$  restricts to an isomorphism $\whvp_{\fG} \colon \fD_{\vp} \to \fD_{\vp}$.
\end{cor}
\proof
We have
 $\ds \whvp_{\fG}(\fD_{\vp}) = \cap_{\ell\geq 0} \ \whvp_{\fG}( \fG_{\ell}) = \cap_{\ell\geq 0} \   \fG_{\ell +1} = \fD_{\vp}$, 
so the restriction $\whvp_{\fG} \colon \fD_{\vp} \to \fD_{\vp}$ is onto.  As $\whvp_{\fG}$ is an embedding, its restriction to $\fD_{\vp}$  is injective.
\endproof

\begin{remark}\label{rmk-strategy}
{\rm
We want to obtain a version of Proposition~\ref{prop-profiniteG} for the completion $\whGamma_{\vp}$ of $\G$, in place of the universal profinite completion $\fG$. That is, we will show that $\vp \colon \G \hookrightarrow \G$ induces an open embedding $\whvp \colon \whGamma_\vp \hookrightarrow \whGamma_\vp$ (see Proposition \ref{prop-embedding} in the next section). 

  Unfortunately, the argument in Proposition \ref{prop-profiniteG} does not directly generalize to the case of the closure of the action $\whGamma_{\vp}$.  Indeed, in the above proof, the key idea is that the system $\cN_1'$ is cofinal in $\cN_1$, and this follows because $\cN_1'$ contains all normal cores (in $\Gamma$) of members of $\cN_1$. 

On the other hand, the group $\whGamma_{\vp}$ is defined as the closure of the action group $\Phi(\G) \subset \Homeo(X_{\vp})$, and the map $\vp$ induces a natural isomorphism with the closure  of the image $\vp(\G) = \G_1 = \G_{U_1}$ in $\Homeo(U_1)$. In order to obtain an embedding of $\whGamma_{\vp}$ into itself, we must relate the closure  of $\G_1$ in $\Homeo(U_1)$ with that in $\Homeo(X_{\vp})$. Thus,  the above algebraic argument using normal cores needs to be replaced by a dynamical argument. The key point is that one needs to show that if $g\in\whGamma_\vp$ acts trivially on $U_1$, then it acts trivially on $X_\vp$. This dynamical regularity, i.e. that the action of $g\in\whGamma_\vp$ is determined by its behavior on any open set, is the goal of the next two sections.
}
\end{remark}

\section{Regularity of Cantor actions}\label{sec-regularity}

 In this section we recall the notion of    \emph{quasi-analytic} actions, and the localized version of this property. 
 This is a type of regularity property for Cantor actions,   introduced in the  works  by   {\'A}lvarez L{\'o}pez,   Candel, and Moreira Galicia  \cite{ALC2009,ALM2016}, inspired by work of Haefliger \cite{Haefliger1985}. 
 We then consider this property for the action $\whvp_{\fG} \colon \fG \times X_{\vp} \to X_{\vp}$  of the profinite completion $\fG$. 
In the following definitions, $H$ denotes a topological group which need not be countable.

  \begin{defn}\label{def-qa}
 An action $\Phi \colon H \times \fX \to \fX$, where 
    $H$ is a topological group and  $\fX$ is a Cantor space,    is said to be \emph{quasi-analytic (or QA)} if for each clopen set $U \subset \fX$, 
if  the action of $g \in H$ satisfies $\Phi(g)(U) = U$ and the restriction $\Phi(g) | U$ is the identity map on $U$, 
  then $\Phi(g)$ acts as the identity on all of $\fX$.  
  \end{defn}
  Note that if an action is not quasi-analytic, then there is some  $g \in H$ whose action $\Phi(g)$ on $\fX$ is non-trivial, yet there is a clopen subset $U$ such that the restriction of the action $\Phi(g)$  to $U$ is the identity, and thus the isotropy group for the action at a point $x \in U$ is non-trivial. So, for example, if  the space  $\fX$ is homeomorphic to a profinite group $\fH$ for which the action $\Phi$ is defined by group multiplication, so that the action is induced by a group homomorphism $\Phi \colon H \to \fH$, then the action is  quasi-analytic. 
A topologically free action, as in Definition~\ref{def-topfree}, is quasi-analytic. Conversely,  the Baire Category Theorem implies that an \emph{effective}  quasi-analytic action of a \emph{countable} group  is topologically free \cite[Section~3]{Renault2008}.

  A local formulation of the QA condition  for   actions was introduced  in the works \cite{DHL2016c,HL2018a}, and has proved very useful for the study of the  dynamical properties of    equicontinuous Cantor actions.
      \begin{defn}  \label{def-LQA}  
   An action $\Phi \colon H \times \fX \to \fX$, where 
    $H$ is a topological group and  $\fX$ a Cantor metric space with metric $\dX$,   is   \emph{locally quasi-analytic} (or \emph{LQA}), if there exists $\e > 0$ such that for any non-empty open set $U \subset \fX$ with $\diam (U) < \e$,  and  for any non-empty open subset $V \subset U$,  if the action of $g \in H$ satisfies $\Phi(g)(V) = V$ and the restriction $\Phi(g) | V$ is the identity map on $V$,    then $\Phi(g)$ acts as the identity on all of $U$.  
\end{defn}
 Examples of minimal equicontinuous Cantor actions   which are locally quasi-analytic, but not quasi-analytic,  are elementary to construct; some examples are given in \cite{DHL2016c,HL2018a}. 
 
If $(\fX, H,  \Phi)$  is an equicontinuous  Cantor action which is not quasi-analytic, then the isotropy group defined in \eqref{iso-defn2}     is non-trivial. On the other hand, there are actions with non-trivial isotropy group   that are quasi-analytic (see Section~\ref{subsec-nilpotent} below, and the examples in \cite{DHL2016c}).
   Finally, we define:

   \begin{defn}  \label{def-stable}  
 An equicontinuous  Cantor action  $(\fX, \G, \Phi)$ is said to be \emph{stable} if the associated profinite action  $\whPhi \colon \overline{\Phi(\G)}  \times \fX \to \fX$ is locally quasi-analytic. The action  is said to be \emph{wild} otherwise. 
\end{defn}
There are many examples of wild Cantor actions. For example, the actions of branch groups on the boundaries of their associated trees are always wild  \cite{GL2019}.  The work \cite{ALBLLN2020} gives the construction of wild Cantor actions exhibiting a variety of characteristic properties, using algebraic methods.

 Here is our main technical result for the profinite actions associated to renormalizations.  
  \begin{thm}\label{thm-proaction}
  Let $\G$ be a finitely generated group and $\vp \colon \G \to \G$ a renormalization of $\G$. Then the profinite action  $\whPhi_{\fG} \colon \fG \times X_{\vp} \to X_{\vp}$ is quasi-analytic.   
  \end{thm}
\proof

Let $g \in \fG$ be such that $\whPhi_{\fG}(g)$ acts non-trivially on $X_{\vp}$.  Suppose there exists a non-empty open set $U \subset X_{\vp}$ such that   $\whPhi_{\fG}(g)$ acts on $U$ as the identity. 

The orbit of every point of $X_{\vp}$ is dense in $X_{\vp}$ under the action of $\G$, so there exists $h \in \G$ such that $\Phi_{\vp}(h)(x_{\vp}) \in U$. Set $g' = h^{-1} g h$ so that $\whPhi_{\fG}(g')$ fixes the open set $U' = \Phi_{\vp}(h^{-1})(U)$. In particular, $\whPhi_{\fG}(g')$ fixes $x_{\vp}$ and hence $g' \in \fD_{\vp}$. Thus, we can assume without loss of generality that   $\whPhi_{\fG}(g)$ acts as the identity on $U$ and $x_{\vp} \in U$, so that $g \in \fD_{\vp}$.

The nested clopen sets $\cU = \{U_{\ell}    \mid \ell \geq 0\}$ form a neighborhood basis at $x_{\vp}$ so there exists some $k_0 > 0$ such that $U_k \subset U$ for all $k \geq k_0$. Thus, for all $k \geq k_0$ we have $\lambda_{\vp}^k(X_{\vp}) \subset U$.

Recall that the restriction $\whvp_{\fG} \colon \fD_{\vp} \to \fD_{\vp}$ is an automorphism  by Corollary~\ref{cor-automorphism}. 
Thus for  $g \in \fD_{\vp}$  there is a well-defined element $g_{\ell} = \whvp_{\fG}^{\ell}(g) \in \fD_{\vp}$   for all $\ell \in \mZ$.

\begin{lemma}\label{lem-key1}
Let $g \in \fD_{\vp}$, and suppose $g$ acts trivially on $U_{k_0}$, for some $k_0 \geq 0$. Then for all $\ell \geq k_0$,  the action of $g_{-\ell} = \whvp_{\fG}^{-\ell}(g) \in \fD_{\vp}$ on $X_{\vp}$ is trivial. % In particular, for $\ell \geq k_0$, $g_{-\ell} \in G_\vp$.
\end{lemma}
\proof

For $x \in X_{\vp}$ and $\ell \geq 0$, set $x_{\ell} = \lambda_{\vp}^{\ell}(x)$. 
Choose $g_x \in \fG$ such that $x = g_x   \fD_{\vp}$ via the identification $X_{\vp} \cong \fG/\fD_{\vp}$;  
that is, $x$ is represented in $X_{\vp}$ by the coset $g_x \fD_{\vp}$. 
Recall that under this identification, for $h \in \fG$ the action of $\whPhi_{\fG}(h)$ on $X_{\vp}$ becomes left multiplication by $h$. That is,
 $\whPhi_{\fG}(h)(x) = h \cdot g_x \fD_{\vp} = h g_x \cdot \fD_{\vp}$. Then for $\ell \geq k_0$ we  have that 
 \begin{equation}\label{eq-calculate1}
x_{\ell} = \lambda_{\vp}^{\ell}(x) = \lambda_{\vp}^{\ell}(g_x \fD_{\vp}) = \whvp_{\fG}^{\ell}(g_x) \fD_{\vp} \in U_{k_0}\ .
\end{equation}
 Thus, for  $\ell \geq k_0$ we have
 \begin{equation}\label{eq-calculate2}
x_{\ell} = \whPhi_{\fG}(g)(x_{\ell}) = \whPhi_{\fG}(g)( \whvp_{\fG}^{\ell}(g_x) \fD_{\vp}) = g \ \whvp_{\fG}^{\ell}(g_x) \fD_{\vp} \ .
\end{equation}
 So for $g \in \fD_{\vp}$ and $\ell \geq k_0$, using that $\whvp_{\fG} \colon \fG \to \fG$ is a homomorphism, we have
 \begin{equation}\label{eq-calculate3}
 x_{\ell} = g \ \whvp_{\fG}^{\ell}(g_x) \fD_{\vp}  = \whvp_{\fG}^{\ell}(\whvp_{\fG}^{-\ell}(g))  \ \whvp_{\fG}^{\ell}(g_x) \fD_{\vp} = \whvp_{\fG}^{\ell}(\whvp_{\fG}^{-\ell}(g) g_x) \fD_{\vp} = \whvp_{\fG}^{\ell}(g_{-\ell} g_x) \fD_{\vp}\ .
 \end{equation}
Thus for $g \in \fD_{\vp}$, $x \in X_{\vp}$  and $\ell \geq k_0$,
 \begin{eqnarray}
x = \lambda_{\vp}^{-\ell}(x_{\ell}) & = &   \lambda_{\vp}^{-\ell}( \whvp_{\fG}^{\ell}(g_{-\ell} g_x) \fD_{\vp})   \label{eq-calculate}\\
& = &  \whvp_{\fG}^{-\ell}(\whvp_{\fG}^{\ell}(g_{-\ell} g_x)) \fD_{\vp}   \nonumber \\
& = &   g_{-\ell} g_x  \fD_{\vp}    =  \whPhi_{\fG}(g_{-\ell})(x) \ . \nonumber
\end{eqnarray}
 That is, $\whPhi_{\fG}(g_{-\ell})(x) = x$ for all $x \in X_{\vp}$ and $\ell \geq k_0$, as was to be shown.
\endproof

 Note that  for $g \in \fG$ the equicontinuous action of $\whvp_{\fG}(g)$ on $X_{\vp}$ is approximated by the action on the finite quotient spaces $X_{\ell}$ for $\ell > 0$. Thus, the assumption that $\whvp_{\fG}(g)$ acts non-trivially on $X_{\vp}$ implies there exists some $m_0 > 0$ such that the induced action of $\whvp_{\fG}(g)$ on $X_{m_0} = \G/\G_{m_0}$ is non-trivial for some $m_0 > 0$. Denote this action  by 
$\whPhi_{m_0}(g) \in \Perm(X_{m_0})$, where $\Perm(X_{m_0})$ is the   group of permutations of the finite set $X_{m_0}$, hence $\Perm(X_{m_0})$ is a finite group. 

The second key observation required for the proof of Theorem~\ref{thm-proaction} is the following ``periodicity'' of the restricted action $\whPhi_{m_0} \colon \fD_{\vp} \to \Perm(X_{m_0})$, which   allows us to promote properties of the action of $g_{-\ell}$ on small scales $U_k$ (for $k$ large) to global properties of the action of $g$.

\begin{lemma}\label{lem-key2}
Let $g \in \fD_{\vp}$. Then 
for every $m_0\geq 1$, there exists $N_0\geq 1$ such that for all $\ell \geq 1$ we have   $\whPhi_{m_0}(g_{-N_0 \ell}) =\whPhi_{m_0}(g)$.
\end{lemma}
\proof
Let $m_0\geq 1$. We note the following two elementary properties of $\text{Hom}(\fG,\Perm(X_{m_0}))$:
\begin{enumerate}
\item For any $\sigma\in\text{Hom}(\fG,\Perm(X_{m_0}))$ and $\ell\geq 0$, we have $\sigma\circ\whvp_{\fG}^\ell\in\text{Hom}(\fG,\Perm(X_{m_0}))$. 
\item $\text{Hom}(\fG,\Perm(X_{m_0}))$ is a finite set, as $\Gamma$ is finitely generated and $\Perm(X_{m_0})$ is   finite.  Indeed, every element in $\text{Hom}(\fG,\Perm(X_{m_0}))$ is a group homomorphism, and so it is determined by its values on the generating set of a dense subgroup $\G$ of $\fG$, which is finite. 
\end{enumerate}
 Now consider $\whPhi_{m_0} \in \text{Hom}(\fG,\Perm(X_{m_0}))$. By properties (1) and  (2) above, 
 there exist $0\leq N_1<N_2$ such that $\whPhi_{m_0}\circ \whvp_{\fG}^{N_2}=\whPhi_{m_0}\circ \whvp_{\fG}^{N_1} \in \text{Hom}(\fG,\Perm(X_{m_0}))$. In particular, their restrictions satisfy  $\whPhi_{m_0}\circ \whvp_{\fG}^{N_2}=\whPhi_{m_0}\circ \whvp_{\fG}^{N_1} \colon \fD_{\vp} \to \Perm(X_{m_0})$.
 
 Now recall from Corollary~\ref{cor-automorphism} that 
   $\whvp_{\fG}$ restricts to an automorphism on $\fD_{\vp}$, so that $\whvp_{\fG}^{-N_1} \colon \fD_{\vp} \to \fD_{\vp}$ is well-defined. It follows that 
   $$\whPhi_{m_0} \circ \whvp_{\fG}^{N_2-N_1} =
   \whPhi_{m_0}\circ \whvp_{\fG}^{N_2} \circ \whvp_{\fG}^{-N_1} = 
  \whPhi_{m_0}\circ \whvp_{\fG}^{N_1} \circ \whvp_{\fG}^{-N_1} = 
   \whPhi_{m_0} \colon \fD_{\vp} \to \Perm(X_{m_0}) \ .$$
    Therefore $N_0=N_2-N_1$ satisfies the conclusion of the lemma.
\endproof

We can now complete the proof that  the action $\whPhi_{\fG} \colon \fG \times X_{\vp} \to X_{\vp}$ is quasi-analytic. If not, then  there  exists  $g \in \fD_{\vp}$   such that $\whPhi_{\fG}(g)$ acts non-trivially on $X_{\vp}$, and    a non-empty open set $U \subset X_{\vp}$ with $x_{\vp} \in U$ such that   $\whvp_{\fG}(g)$ acts on $U$ as the identity. 

Let  $m_0 > 0$ be such that the induced action of $\whvp_{\fG}(g)$ on $X_{m_0} = \G/\G_{m_0}$ is non-trivial. 
 
Then there exists some $k_0 > 0$ such that $U_k \subset U$ for all $k \geq k_0$. We assume that $k_0 \geq m_0$.

By Lemma~\ref{lem-key1}, for all $\ell \geq k_0$ we have $\whPhi_{m_0}(g_{-\ell})(x) = x$ for all $x \in X_{m_0}$.

By Lemma~\ref{lem-key2}, there exists $N_0 > 0$ so that $\whPhi_{m_0}(g_{-N_0 \ell}) =\whPhi_{m_0}(g)$ for all $\ell \geq 0$. 

However, for $N_0 \ell > k_0$ we obtain a contradiction, as $\whPhi_{m_0}(g)$ is assumed to act non-trivially on $X_{m_0}$ while $\whPhi_{m_0}(g_{-N_0 \ell})$ acts trivially on $X_{m_0}$.  

Thus,  the action of $\whPhi_{\fG}$ must be quasi-analytic.
   \endproof

Finally,  note that   Theorem~\ref{thm-proaction} shows that the profinite action   $(X_{\vp}, \fG, \whPhi_{\fG})$ is quasi-analytic, so the same holds for the   action $(X_{\vp}, \G, \Phi_{\vp})$ obtained by restricting  the action to the image   $\Phi_\vp(\G) \subset \fG$. Then the 
 Baire Category Theorem implies (see  \cite[Proposition~2.2]{HL2018b} for example) that   if the action $\Phi_{\vp}$ is effective, that is if $\Phi_{\vp} \colon \G \to \Homeo(X_{\vp})$ is injective,  then the   action $(X_{\vp}, \G, \Phi_{\vp})$ is topologically free, as asserted in Theorem~\ref{thm-main1}.

  \begin{cor}\label{cor-stable}
  Let $\G$ be a finitely generated group and $\vp \colon \G \to \G$ a renormalization of $\G$. Then the profinite action  $\whPhi_{\vp} \colon \whGamma_{\vp} \times X_{\vp} \to X_{\vp}$ is quasi-analytic, and the action 
 $(X_{\vp}, \G, \Phi_{\vp})$ is  stable.   
\end{cor}
\proof
Let   $U \subset X_{\vp}$ be a clopen set and $\whg \in \whGamma_{\vp}$ such that $\whPhi_{\vp}(\whg)$ restricts to the identity on $U$. 
Recall that $\Pi_{\vp} \colon \fG \to \whGamma_{\vp}$  is onto, so there exists  $g \in \fG$   such that $\Pi_{\vp}(g) = \whg$.
Then the action $\whPhi_{\fG}(g)$ restricts to the identity on $U$, so $\whPhi_{\fG}(g)$ acts as the identity on $X_{\vp}$ by Theorem~\ref{thm-proaction}. % Thus, $\whg$ is the identity, as required to show.
\endproof

\section{Open embeddings}\label{sec-open}
  
 In this section, given a renormalization $\vp \colon \G \to \G$ with associated Cantor action $(X_{\vp}, \G, \Phi_{\vp})$, 
 we obtain  a structure theory for the profinite group $\whGamma_{\vp}$ that is the key to the proof of Theorem~\ref{thm-main4}.
Recall that  Proposition~\ref{prop-profiniteG} showed that the induced map $\whvp_{\fG} \colon \fG \to \fG$ of the profinite completion $\fG$ of $\G$ is an open embedding. We thus obtain by  \cite[Theorem~3.10]{WvL2018b}, which is a reformulation of the results of Reid in   \cite{Reid2014}, the following structure theorem:

 \begin{thm} \label{thm-reid}There exist closed subgroups $C_{\vp} \subset \fG$ and $Q_{\vp} \subset \fG$  so that:
\begin{enumerate}
\item $\fG \cong C_{\vp} \rtimes Q_{\vp}$, where $C_{\vp}$ is a pro-nilpotent normal subgroup of $\fG$;
\item $C_{\vp}$ is $\whvp_{\fG}$-invariant, and $\whvp_{\fG}$ restricts to an open contracting embedding on $C_{\vp}$;
\item $Q_{\vp}$ is $\whvp_{\fG}$-invariant, and $\whvp_{\fG}$ restricts to an automorphism of $Q_{\vp}$.
\end{enumerate}
Moreover, let $\whe \in \fG$ be the identity element, then we have
\begin{equation}\label{eq-characterizations}
C_{\vp} = \{g \in  \fG  \mid \lim_{n \to \infty} \ \whvp_{\fG}^n(g) = \whe\} \quad , \quad Q_{\vp} = \bigcap_{n > 0} \ \whvp_{\fG}^n(\fG) ~ .
\end{equation}
\end{thm}
  
Next, we show   that $\whvp_\fG$  induces an open embedding   $\whvp \colon \whGamma_{\vp} \to \whGamma_{\vp}$ as promised in Remark~\ref{rmk-strategy}.
  \begin{prop}\label{prop-embedding}
    Let $\vp$ be a renormalization of the finitely-generated group $\G$. Then $\vp$ induces an injective  homomorphism  $\whvp \colon \whGamma_{\vp} \to \whGamma_{\vp}$ whose image is a clopen subgroup of $\whGamma_{\vp}$.
  \end{prop}
  \proof

Let   $\fN(\Pi_{\vp}) = \ker \{ \Pi_{\vp} \colon \fG \to \whGamma_{\vp}\}   \subset \fG$ be  the kernel of   the homomorphism   $\ds \Pi_{\vp}$.   
We claim  that $\whvp_{\fG} \colon \fG \to \fG$ descends to a homomorphism 
  \begin{equation}\label{eq-inducedmap}
\whvp \colon \whGamma_{\vp} \cong \fG/\fN(\Pi_{\vp}) \longrightarrow \whGamma_{\vp} \cong \fG/\fN(\Pi_{\vp}) \ .
\end{equation}
   
   Recall that 
Theorem~\ref{thm-main1} implies  that for   $g \in \fG$,  the action of $\whvp_{\fG}(g)$ on $X_{\vp}$ 
  is locally determined.  In particular, the action $\whvp_{\fG}(g)$  is determined by its restriction to the clopen subset $U_1 = \lambda_{\vp}(X_{\vp})$.

  For $g \in \fN(\Pi_{\vp})$, observe that $\whvp_{\fG}(g)$ acts as the identity on the clopen subset $U_1 = \lambda_{\vp}(X_{\vp})$. As the action $\whPhi_{\fG}$ is quasi-analytic, 
  this implies that $\whvp_{\fG}(g)$ acts as the identity on  $X_{\vp}$, and thus $\whvp_{\fG}(g) \in \fN(\Pi_{\vp})$. 
That is,  $\whvp_{\fG}(\fN(\Pi_{\vp})) \subset \fN(\Pi_{\vp}) \subset \fG$, and thus we have the composition of homomorphisms
  \begin{equation}
 \whvp \colon  \whGamma_{\vp} = \fG/\fN(\Pi_{\vp}) \to \whvp(\fG)/\whvp(\fN(\Pi_{\vp})) \to \fG/\whvp(\fN(\Pi_{\vp}) \to \fG/\fN(\Pi_{\vp}) = \whGamma_{\vp} \ 
\end{equation}
which defines the    map  \eqref{eq-inducedmap}.
We claim that  $\whvp$ is injective. If not, let $\gamma \in \whGamma_{\vp}$ such that $\whvp(\gamma) = {\rm Id}$.
That is, $\whvp(\gamma) \in \whGamma_{\vp}$ acts as the identity on $X_{\vp}$. In particular, $\whvp(\gamma)$ acts as the identity on $\lambda_{\vp}(X_{\vp})$, so for $x \in X_{\vp}$, 
$$\lambda_{\vp}(x) = \whvp(\gamma) \cdot \lambda_{\vp}(x) = \lambda_{\vp}(\gamma x) \ .$$
    As $\lambda_{\vp}$ is an injection, we have $\gamma  x = x$ for all $x \in X_{\vp}$, and thus $\gamma \in \Homeo(X_{\vp})$ is the identity, as was to be shown.
   \endproof

We use  the conclusions of Theorem~\ref{thm-reid} and Proposition~\ref{prop-embedding} to obtain:  
 
 \begin{thm} \label{thm-main2}
 Let $\vp \colon \G \to \G$ be a renormalization for the  finitely generated group $\G$, with associated Cantor action $(X_{\vp}, \G, \Phi_{\vp})$.      Let  $\whvp \colon \whGamma_{\vp} \to \whGamma_{\vp}$ be the   embedding induced from $\vp$. Then
  there exists a closed \emph{pro-nilpotent normal} subgroup  $\whN_{\vp} \subset \whGamma_{\vp}$ so that:
\begin{enumerate}
\item $\whGamma_{\vp} \cong \whN_{\vp} \rtimes \cD_{\vp}$ is a semi-direct product;
\item $\whN_{\vp}$ is $\whvp$-invariant, and $\whvp$ restricts to an open contracting embedding on $\whN_{\vp}$;
\item $\cD_{\vp}$ is $\whvp$-invariant, and $\whvp$ restricts to an automorphism of $\cD_{\vp}$.
\end{enumerate}
Moreover, let $\whe \in \whGamma_{\vp}$ be the identity element, then we have
\begin{equation}\label{eq-contractingdef}
\whN_{\vp} = \{g \in  \whGamma_{\vp}  \mid \lim_{n \to \infty} \ \whvp^n(g) = \whe\} \quad , \quad \cD_{\vp} = \bigcap_{n > 0} \ \whvp^n(\whGamma_{\vp}) ~ .
\end{equation}
\end{thm}
 \proof
Recall that by Theorem~\ref{thm-reid} the embedding   $\whvp_{\fG} \colon \fG \to \fG$ induces an isomorphism   
  $\fG \cong C_{\vp} \rtimes Q_{\vp}$, where $C_{\vp}$ and $Q_{\vp}$ are characterized by the formulae in \eqref{eq-characterizations}.
   First, we show:
\begin{lemma}\label{lem-disc}
$\cD_{\vp} =   \Pi_{\fG}(Q_{\vp}) \subset \whGamma_{\vp}$.
\end{lemma}
\proof
Recall that  the clopen neighborhoods $U_{\ell}$ of $x_{\vp}$ are defined by \eqref{eq-clopenbasis}, and   for each $\ell \geq 0$ we have   $U_{\ell} = \lambda_{\vp}^{\ell}(X_{\vp})$.
For each $\ell \geq 0$, define the clopen subset  $\whU_{\ell} = \{ \gamma \in \whGamma_{\vp} \mid \whPhi_{\vp}(\gamma)(U_{\ell}) = U_{\ell} \} \subset \whGamma_{\vp}$.

Also, recall that $\cD_{\vp} = \{\gamma \in \whGamma_{\vp} \mid \gamma \cdot x_{\vp} = x_{\vp}\}$.   
As $x_{\vp} =  \bigcap_{\ell \geq 0} \ U_{\ell}$,    we then have $\cD_{\vp} = \bigcap_{\ell \geq 0} \ \whU_{\ell}$, and so $\whU_{\ell} = \whvp^{\ell}(\whGamma_{\vp})$ where 
  $\whvp \colon \whGamma_{\vp} \to \whGamma_{\vp}$ was defined in  Proposition~\ref{prop-embedding}.

Recall that $\fG_{\ell} = \whvp_{\fG}^{\ell}(\fG) \subset \fG$,  and thus $\whU_{\ell} = \Pi_{\fG}(\fG_{\ell})$. Then we have
\begin{equation}\label{eq-dcalc}
\Pi_{\fG}(Q_{\vp}) = 
 \Pi_{\fG} \left\{ \bigcap_{\ell \geq 0} \ \fG_{\ell} \right\}  = 
  \bigcap_{\ell \geq 0} \ \Pi_{\fG}(\fG_{\ell}) = 
 \bigcap_{\ell \geq 0} \ \whU_{\ell} = \cD_{\vp} \ ,
 \end{equation}
as was to be shown.
\endproof

  Next, set   $\whN_{\vp} = \Pi_{\fG}(C_{\vp}) \subset \whGamma_{\vp}$ which is a pro-nilpotent closed subgroup. Then by an argument exactly analogous to the proof of Lemma~\ref{lem-disc}, we have 
\begin{equation}
\whN_{\vp} = \{\gamma \in \whGamma_{\vp} \mid \lim_{\ell \to \infty} \ \whvp^{\ell}(\gamma) = \whe \in \whGamma_{\vp}\} \ .\end{equation}
The proof of Theorem~\ref{thm-main2} now follows.
\endproof
  
 Note that the identities \eqref{eq-contractingdef} in Theorem~\ref{thm-main2} identify the images of the groups $C_{\vp}$ and $Q_{\vp}$ in $\Homeo(X_{\vp})$ in terms of the dynamical properties of the action $\whvp$ on $\whGamma_{\vp}$.

 The conclusions of Theorem~\ref{thm-main2} are illustrated in various examples of renormalizable groups and self-embeddings in     Section~\ref{sec-examples}, and also by the examples in the works \cite{NekkyPete2011,WvL2018a,WvL2018b}.  Moreover, the conclusion    that  $\vp$ induces an automorphism of the discriminant group $\cD_{\vp}$     has applications to the constructions of examples of Cantor actions using the \emph{Lenstra method} as given in \cite[Section~8.2]{HL2018a}.

\section{Finite discriminant}\label{sec-finitedisc}

 We next consider the consequences of  Theorem~\ref{thm-main2}, for  $\cD_{\vp}$   a finite group  and prove Theorem \ref{thm-main4}. 
 
 \proof \emph{(of Theorem \ref{thm-main4})}
 We first assume that the discriminant group $\cD_{\vp}$   is trivial, and show that  the quotient group $\G/C(\cG_{\vp})$ is nilpotent, where $C(\cG_{\vp})$ is  the normal core of the intersection  $K(\cG_{\vp}) = \bigcap_{\ell \geq 0} \G_\ell \subset \G$, as defined in \eqref{eq-kernels}.    Recall that $C(\cG_{\vp}) \subset \G$ is identified with the kernel of the  homomorphism 
  $\Phi_{\vp} \colon \G \to \whGamma_{\vp} \subset \Homeo(X_{\vp})$, and that    $\whvp \colon \whGamma_{\vp} \to \whGamma_{\vp}$ was defined in  Proposition~\ref{prop-embedding}.

 Note that $\vp$ restricts to an isomorphism of $K(\cG_{\vp})$ by its definition, and so $\vp$  also maps $C(\cG_{\vp})$ isomorphically to itself, and thus induces an embedding $\vp' \colon \G/C(\cG_{\vp}) \to \G/C(\cG_{\vp})$.
  Then without loss of generality, we can replace $\G$ with $\G/C(\cG_{\vp})$, so can assume that $\Phi_{\vp} \colon \G \to \whGamma_{\vp}$ is an embedding, and identify $\G$ with its image $\Phi_{\vp}(\G)$. 
As we assume that  $\cD_{\vp}$   is trivial,  by Theorem~\ref{thm-main2} we have 
  $\whGamma_{\vp} \cong \whN_{\vp}$ where $\whN_{\vp}$ is  a closed pro-nilpotent normal group.
  
 Section 3 of the work \cite{WvL2018b} gives an overview of  some of the structure theory of pro-nilpotent groups, and we recall those aspects as required for the proof of Theorem~\ref{thm-main4}. 
   First, $\whN_{\vp}$ admits a splitting by \cite[Theorem B]{GW2010} as $\whN_{\vp} \cong   \whN_{\infty} \times \whN_{tor}$ where $\whN_{\infty}$ is a torsion-free  nilpotent group  and $\whN_{tor}$ is a torsion group with   bounded exponent, by results of \cite{KW1993}.   We now claim:
\begin{lemma}\label{lem-trivialtor}
If $\cD_\vp$ is trivial, then $\whN_{tor}$ is the trivial group.
\end{lemma}
\proof
Let $\pi_{tor} \colon \whN_{\vp} \to \whN_{tor}$ be the projection, then the image $\pi_{tor}(\G) \subset \whN_{tor}$ is dense.

The abelianization $\whA_{tor}$ of $\whN_{tor}$  is an abelian group of bounded exponent, which is trivial if and only if  $\whN_{tor}$ is trivial.  By Pr\"ufer's First Theorem (see \S~24 of \cite{Kurosh1960}),   $\whA_{tor}$ is a direct sum  of (possibly infinitely many) cyclic groups.  
As $\G$  is finitely generated,  the image of  $\G$ in       $\whA_{tor}$     is finite rank  and dense, and therefore the abelianization    $\whA_{tor}$ has finite rank. Thus, $\whA_{tor}$ is a direct sum of finitely many cyclic groups, hence  is a  finite group. 
 
 Note that the contraction mapping $\whvp \colon \whN_{\vp} \to \whN_{\vp}$ induces a contraction mapping 
  $\whvp_{tor} \colon \whN_{tor} \to \whN_{tor}$.
  
 The second part of Theorem~B in  Gl\"ockner and Willis \cite{GW2010}   proves the existence of a Jordan-H\"older series for bounded exponent contraction groups with each composition factor a \emph{simple} contraction group. Here we say a contraction group with contraction $\alpha$ is \emph{simple} if it has no nontrivial, proper closed normal $\alpha$-invariant subgroup. Further, the simple contraction groups are classified as shifts on $F^{\mathbb N}$ where $F$ is a finite simple group. By considering the first composition factor, we see that $\whN_{tor}$ has a quotient of the form $F^{\mathbb N}$ where $F$ is a finite simple group. Since $\whN_{tor}$ is solvable of bounded exponent \cite{Reid2014}, we conclude    that $F$ is abelian. In particular $\whN_{tor}$ has an infinite abelian quotient, which contradicts the fact that $\whA_{tor}$ is  a  finite group, as shown previously.
 It follows that $\whN_{tor}$ must be  the trivial group
 \endproof

    Since by assumption $\cD_\vp$ is trivial, then $\G \subset  \whN_{\infty}$. Now observe that by Lemma~\ref{lem-trivialtor}, $\whN_{\infty}$  is a torsion-free  nilpotent group, thus $\G$ is nilpotent. This concludes the proof of Theorem~\ref{thm-main4} in the case where $\cD_{\vp}$ is trivial.

 Next,  assume that $\cD_{\vp}$ is a finite group.
 By Theorem~\ref{thm-main2}, we have  $\cD_{\vp} =   \Pi_{\fG}(Q_{\vp}) \subset \whGamma_{\vp}$ and its intersection with  $\whN_{\vp} = \Pi_{\fG}(C_{\vp})$ is the trivial subgroup. As $\cD_{\vp}$ is a finite group, it follows that $\whN_{\vp}$ is a clopen subset of $\whGamma_{\vp}$, and so   $\Lambda_{\vp} \equiv \G \cap \whN_{\vp}$   is a dense subgroup of $\whN_{\vp}$

The restriction of  $\whvp$ defines   a contraction mapping $\whvp \colon \whN_{\vp} \to \whN_{\vp}$. We can thus   apply the above arguments for     the case when $\cD_{\vp}$ is trivial to the action of $\Lambda_{\vp}$ on $\whN_{\vp}$ to conclude 
  that  $\Lambda_{\vp}$ is nilpotent. As $\Lambda_{\vp}$ has     finite index in $\G$,   this 
  completes   the proof of Theorem~\ref{thm-main4}.

 \proof \emph{(of Theorem~\ref{thm-contractionprinciple})} Assume that both 
  $\cD_{\vp}$  and the subgroup $K(\cG_{\vp})$ in \eqref{eq-kernels} are finite groups.   Thus its core $C(\cG_{\vp}) \subset K(\cG_{\vp})$ is also finite.
  Recall that in the above proof of Theorem~\ref{thm-main4}, we replaced $\G$ with the quotient $\G/C(\cG_{\vp})$, and concluded that   $\G/C(\cG_{\vp})$ contains a nilpotent subgroup of finite index.    
  In the case where both groups $\cD_{\vp}$ and $K(\cG_{\vp})$ are trivial, then the claim of the corollary follows directly from  Theorem~\ref{thm-main4} and Theorem~\ref{thm-main2}.
  In the case where both groups are finite, we have that 
 $C(\cG_{\vp})$ is a finite normal subgroup of $\G$ and $\G/C(\cG_{\vp})$ contains a nilpotent subgroup of finite index, which implies that     $\G$   contains a nilpotent subgroup of finite index.
  This completes   the proof. \endproof

  \begin{remark}
  {\rm
  We clarify the relation between the groups   $K(\cG_\vp)$ and $\cD_{\vp}$  in the hypothesis of Theorem~\ref{thm-contractionprinciple}. The group $K(\cG_\vp)$ contains every element of $\G$ which fixes the base point $x_\vp$, while the group $\cD_\vp$ contains every element in $\whGamma_\vp$ which fixes $x_\vp$. The relationship between $K(\cG_\vp)$ and $\cD_\vp$ is provided by an embedding $\G \to \whGamma_\vp \subset \Homeo(X_\vp)$. It follows that if $\cD_\vp$ is a finite group, then the quotient group $K(\cG_\vp)/C(\cG_\vp)$ must be finite, where $C(\cG_\vp)$ is the normal core of $K(\cG_\vp)$. Thus if $\cD_{\vp}$ is finite, $K(\cG_\vp)$ may still be infinite if its normal core $C(\cG_\vp)$ is infinite. Thus both assumptions in \eqref{eq-kernels} in Theorem~\ref{thm-contractionprinciple} are necessary.
  }
  \end{remark}

 \section{Renormalizable Cantor actions}\label{sec-renormalization}

 In this section, we  introduce the notions of   \emph{(virtually) renormalizable}   Cantor actions, and study their regularity properties and invariants, yielding a proof  of Theorem~\ref{thm-dichotomy}.

\subsection{Renormalizable actions} 
 
For a Cantor action   $(\fX,\G,\Phi)$  and an adapted   set $U \subset \fX$, note that    $H_U = \Phi(\G_U) \subset \Homeo(U)$ acts faithfully on $U$, so $(U, H_U, \Phi_U)$ is always an effective action. 
 \begin{defn}\label{defn-renormalizable}
A Cantor action   $(\fX,\G,\Phi)$  is \emph{renormalizable} if it is equicontinuous, and there exists an adapted proper clopen   set $U \subset \fX$  such that   the actions   $(\fX,\G,\Phi)$ and 
$(U, H_U, \Phi_U)$ are isomorphic (as in Definition~\ref{def-isomorphism}) by a homeomorphism $\lambda \colon \fX \to U$ and   isomorphism $\Theta \colon \G \to H_U$, and the intersection $\cap_{\ell \geq 0} \ \lambda^{\ell}(\fX)$ is a point.
\end{defn}

For example, let $(X_{\vp}, \G, \Phi_{\vp})$ be the Cantor action associated to a renormalization  $\vp$ of $\G$.  Suppose the action is effective, then by Theorem \ref{thm-main1} it is topologically free, and thus the   map $\Phi_U \colon \G_U \to H_U$ is an injection. Then  the action is renormalizable with    $\lambda = \lambda_{\vp}$ and $\Theta = \Phi_U \circ \vp \colon \G \to H_U$. 
In general, for a renormalizable action,  there is no requirement that the map $\Phi_U \colon \G_U \to H_U$ is injective, and so   $H_U$ need not be identified with a subgroup of  $\G$.

 \begin{defn}\label{defn-Vrenormalizable}
A Cantor action   $(\fX,\G,\Phi)$  is \emph{virtually renormalizable} if it is equicontinuous, and there exists an adapted     set $V \subset \fX$  such that   the restricted action  $(V, H_V, \Phi_V)$ is renormalizable.
\end{defn}

The class of virtually renormalizable actions is much more general than the renormalizable actions, as it allows for the case when the action map $\Phi \colon \G \to\Homeo(\fX)$ has a non-trivial kernel.
In the following, we show some properties of these actions. We first show:

\begin{prop}\label{prop-qar}
Suppose that  the Cantor action    $(\fX,\G,\Phi)$ is renormalizable   and  locally   quasi-analytic, then it is quasi-analytic. 
\end{prop}
\proof
We assume there is given a homeomorphism $\lambda \colon \fX \to U$ and group isomorphism $\Theta \colon \G \to H_U$ implementing an isomorphism  
of $(\fX,\G,\Phi)$ with $(U,H_U,\Phi_U)$ as in \eqref{eq-isomorphism}.

First, suppose that  the map $\Phi_U \colon \G_U \to H_U$ is injective, and hence is an isomorphism, as it is onto by the definition of $H_U$.
Then  the composition $\vp \equiv \Phi_U^{-1} \circ \Theta \colon \G \to \G$ is a proper inclusion with image $\G_U \subset \G$. As $U$ is adapted, $\G_U$ has finite index in $\G$, and thus  $\vp$ is a renormalization of $\G$. It follows from Theorem \ref{thm-main1} that the action $\Phi$ is quasi-analytic.

It thus suffices to show that if  $\Phi_U \colon \G_U \to H_U$ has a non-trivial kernel $K_U \subset \G$, then   the action $\Phi$ is not locally quasi-analytic, which yields a contradiction. 
We show this using a recursive argument.

Set $U_0 = \fX$, then $U_1  = \lambda(\fX)$ is a clopen set by assumption. Also define clopen sets      $U_{\ell} = \lambda^{\ell}(U_0)$ for $\ell > 1$, so that $U_{\ell} \subset U_{\ell-1}$. 
The assumption in Definition~\ref{defn-renormalizable} that the intersection $\cap_{\ell \geq 0} \ \lambda^{\ell}(\fX)$ is a point, denoted by   $x_{\lambda} \in \fX$,  implies that    $\{U_{\ell}  \mid \ell \geq 0\}$ is an adapted neighborhood basis at $x_{\lambda}$.

Now set $\G_{\ell} = \G_{U_{\ell}}$ for $\ell \geq 0$, and let  $H_{\ell} = \Phi_{U_{\ell}}(\G_{\ell}) \subset \Homeo(U_{\ell})$ for $\ell \geq 0$. Let     $\{\G_{\ell} \mid \ell \geq 0\}$  be the  associated group chain. 
Recall that as the action of $H_1$ on $U_1$ is effective, and the actions $(\fX,\G,\Phi)$ and 
$(U_1, H_1, \Phi_{U_1})$ are isomorphic, so the  action of $\G$ on $\fX$ is effective. That is,   the kernel $K_0 \subset \G$ of $\Phi$ is trivial, and   $\Phi \colon \G \to H_0$ is an isomorphism onto.
  To avoid cumbersome notation, we will identify $\G = H_0$ and write the action as $g \cdot x = \Phi(g)(x)$.

 Now observe that  
$$\G_{\ell+1} = \{g \in \G \mid g \cdot U_{\ell+1} = U_{\ell+1}\} = \{g \in \G_{\ell} \mid g \cdot U_{\ell+1} = U_{\ell+1}\} = (\G_{U_{\ell}})_{U_{\ell+1}}$$
  since $g \cdot U_{\ell+1} = U_{\ell+1}$ implies $g \cdot U_{\ell} = U_{\ell}$, as $U_{\ell}$ is an adapted clopen set and $U_{\ell+1} \subset U_{\ell}$.

  We give the first step of the recursive argument. Define
  \begin{equation}\label{eq-K1}
K_1 ~ \equiv ~ \ker \ \left\{ \Phi_{U_1} \colon \G_1 \to H_1 \subset \Homeo(U_1) \right\} ~ \subset ~ \G_1 \subset \G \ .
\end{equation}
By assumption, the subgroup $K_1$  is non-trivial.

 Let $\Phi^1_{U_1} \colon H_1 \times U_1 \to U_1$ denote the action of $H_1$,  
 and let  $(H_1)_{U_2} \subset H_1$ denote  the elements of $H_1$ which map $U_2$ to itself.
Then introduce the subgroup  $K'_2 \subset (H_1)_{U_2}$ of elements which restrict to  the identity on $U_2$. 
    Then we have:
 \begin{eqnarray}
\quad K'_2 = \ker \left\{\Phi^1_{U_2} \colon (H_1)_{U_2} \to \Homeo(U_2) \right\} 
& = &  \ker \left\{  \Phi^1_{U_2} \colon (H_1)_{\lambda(U_1)} \to \Homeo(\lambda(U_1))   \right\}  \label{eq-recursiveker}\\ 
 & = &  \ker \left\{  \Phi^1_{U_2} \colon \Theta(\Gamma)_{\lambda(U_1)} \to \Homeo(\lambda(U_1))   \right\}   \nonumber
\\
& = &  \Theta  \left( \ker \left\{\Phi_{U_1} \colon \G_{U_1} \to \Homeo(U_1) \right\} \right) = \Theta(K_1)   \ , \nonumber
\end{eqnarray}
where the last equality follows using the isomorphism of   $(\fX,\G,\Phi)$ with  $(U_1,H_1,\Phi_{U_1})$.

By assumption $K_1$ is a non-trivial subgroup, so by \eqref{eq-recursiveker}  we have 
  $K_2' =  \Theta(K_1)$  is also non-trivial. 
  That is,   if  $g \in K_1 \subset   \G_1$ is not the identity, then $g$ acts non-trivially on $U_0 = \fX$ and restricts to the identity on $U_1$ by the definition \eqref{eq-K1} of $K_1$.  
   Thus, $h = \Theta(g) \in H_1$ acts non-trivially on $U_1$ and restricts to the identity on $U_2$.   
   Since $H_1 = \Phi_{U_1}(\G_1)$, there exists $g' \in \Gamma_1$ such that $\Phi_{U_1}(g') = h$. We have found $g' \in \G_1$, such that $g' \notin K_1$ and $g' \in K_2$. Therefore, $K_1$ is a non-trivial proper subgroup of $K_2$. 
     
  Set  $K_{\ell} = \ker \{\Phi_{U_{\ell}} \colon \G_{\ell} \to \Homeo(U_{\ell}) \}$ for $\ell \geq 2$, then by repeating the above arguments in \eqref{eq-recursiveker}, we have  $K_{\ell} \subset K_{\ell+1} \subset \G$ is a proper inclusion for all $\ell \geq 1$. 
    As the diameter of the sets $U_{\ell}$ tends to $0$ as $\ell$ increases, given any adapted set $V \subset \fX$ for the action $\Phi$,  there exist  $\ell > 0$ and $\gamma \in \G$ such that $O = \gamma \cdot U_{\ell} \subset V$. This implies that 
    the dynamics of $\G_{\ell}$ acting on $U_{\ell}$ is conjugate to the restricted action of $\G_V$ on the adapted clopen set $O$.  Thus, 
 there exists some element   $\gamma' \in \G$ such that $\gamma' \cdot O = O$ and the action of $\Phi(\gamma')$ restricted to $O$ is non-trivial, but restricts to the identity on some open set that is a translate of $U_{\ell +1}$ in $O$,  namely $\gamma' = \gamma^{-1} \circ s \circ \gamma$, where $s \in K_{\ell+1}$ and $s \notin K_\ell$. As this holds for all $\ell > 0$,    the action $\Phi$ is not locally quasi-analytic.
 \endproof

We have the following consequence of the above proof of Proposition~\ref{prop-qar}.
 \begin{prop}\label{prop-qar2}
Suppose that  the Cantor action    $(\fX,\G,\Phi)$ is renormalizable   and  locally   quasi-analytic, then  the action is isomorphic to an action $(X_{\vp}, \G, \Phi_{\vp})$ associated to a renormalization $\vp \colon \G \to \G$, and in particular  $\G$ is renormalizable, and  the action $(\fX,\G,\Phi)$ is stable.
\end{prop}

\proof
    
  As in the proof of Proposition~\ref{prop-qar}, let $ \{U_{\ell}  \mid \ell \geq 0\}$ be an adapted neighborhood basis at $x_\vp$, and let  $\{\G_{\ell} \mid \ell \geq 0\}$ is the associated group chain. The action $(\fX,\G,\Phi)$  is quasi-analytic by Proposition~\ref{prop-qar}, so 
we have isomorphisms $\G_{\ell} \cong H_{\ell}$, and in particular 
   the composition $\vp \equiv \Phi_{U_1}^{-1} \circ \Theta \colon \G \to \G$ is a proper inclusion with image $\G_1 \subset \G$ a subgroup of finite index. It then follows that $\G_\ell = \vp^{\ell}(\G)$, and by the results in Section~\ref{sec-renormalizable} (see also \cite{CortezPetite2008,DHL2016a,DHL2016c}) the Cantor action   $(\fX,\G,\Phi)$  is isomorphic to the action $(X_{\vp}, \G, \Phi_{\vp})$. Then the action $(\fX,\G,\Phi)$ is stable by Corollary~\ref{cor-stable}. 
\endproof

As a consequence of the above, we have the following result,   which implies Theorem~\ref{thm-dichotomy}.
 \begin{thm}\label{thm-qar}
Let $\vp$ be a renormalization of $\G$, then   $(X_{\vp}, \G, \Phi_{\vp})$ is virtually renormalizable. Conversely, 
suppose that  a minimal equicontinuous Cantor action    $(\fX,\G,\Phi)$ is renormalizable   and  locally   quasi-analytic, then 
$\G$ is renormalizable, and there is a renormalization $\vp$ such that  $(\fX,\G,\Phi)$ is isomorphic to    $(X_{\vp}, \G, \Phi_{\vp})$.
\end{thm}

\subsection{Renormalizable actions on trees}\label{subsec-renormalizable}
 We will now discuss the relationship between renormalizable actions on Cantor sets and self-similarity properties of groups acting on rooted trees. Recall that, given a group chain $\{\G_\ell\}_{\ell \geq 0}$ consisting of finite index subgroups of $\G = \G_0$, the tree model is a natural action of $\G$ on a rooted tree (see Section~\ref{subsec-treemodels}). We start by briefly recalling this construction. As in the rest of the paper, the action of $\G \subset \Aut(T)$ on the boundary of $T$ is assumed minimal in this section.

Recall that $V_\ell = X_\ell = \G/\G_\ell$, and $V = \bigsqcup_{\ell \geq 0} V_\ell$ is the vertex set of a tree $T$. The boundary $\partial T$ of $T$ consists of all infinite connected paths $(v_\ell)_{\ell \geq 0} \in \partial T$.  Let $\whe = (e_\ell)_{\ell \geq 0} = (e \cdot \G_\ell)$ be the path passing through the coset of the identity $e \in \G$ at each level $V_\ell$.  Then, as in \eqref{eq-clopenbasis}, for $k \geq 0$ the set 
 \begin{align}\label{eq:uk}U_k = \{(w_\ell)_{\ell \geq 0} \in \partial T \mid w_\ell = e \cdot \G_\ell, \, 0 \leq \ell \leq k \} = \{(w_\ell)_{\ell \geq 0} \in \partial T \mid w_k = e \cdot \G_k\} \end{align}
  is a clopen neighborhood of $\whe$ and $\G_k$ is the stabilizer subgroup of $U_k$. Since the action of $\G$ on $\partial T$ is minimal, the induced action of $\G$ on each vertex level $V_\ell$  is transitive, for $\ell \geq 0$.

We now discuss Definition \ref{defn-renormalizable} of a renormalizable action as applied to actions on rooted trees, described in the previous paragraph. Suppose the action $(\partial T, \G)$ is renormalizable with $U = U_1$. The set $U_1$ contains all infinite paths in $\partial T$ which pass through the vertex $e_1 = e\cdot \G_1$, and every element in $\G_1$ fixes $e_1$. For each $g \in \G_1 \subset \Aut(T)$, denote by $g_1 = g|U_1$ the restriction.  If the action $(\partial T,\G)$ is quasi-analytic, then there is precisely one element $g \in \G$ which restricts to $g_1$, and so the map $\Phi_{U_1}: \G_1 \to H_{U_1}$ is invertible. As discussed at the beginning of this section, in this case there is an injective homomorphism $\varphi = \Phi_{U_1}^{-1} \circ \Theta: \G \to \G_1 \subset \G$. In particular, the conditions that $(\partial T, \G)$ is quasi-analytic and minimal implies that there are no elements in $\G$ whose support is contained entirely in $U_1$, and it follows that the class of groups which admit renormalizations does not contain weakly branch groups (see \cite{BGS2012,Grigorchuk2011} for more details about weakly branch groups). 

If the action $(\partial T, \G)$ is renormalizable but not quasi-analytic, then the elements in $H_{U_1}$ can be extended from $U_1$ to $\partial T$ in multiple ways, and $\G$ may have elements with support contained entirely in $U_1$. Such renormalizable actions are wild. Some actions of branch groups belong to this class, for instance, the action of the Grigorchuk group, as we show below.

We will show that renormalizability of the action of $\G$ on $\partial T$ is closely related to the \emph{self-replicating} property of the action of $\G$ on $T$. Given a vertex $v_\ell \in V$, let $v_\ell T$ be a subtree of $T$ with root $v_\ell$. Fix an isomorphism of rooted trees $p_{v_\ell}: T \to v_\ell T$ and let $\partial p_{v_\ell}:\partial T \to \partial (v_\ell T)$ be the induced homeomorphism of boundaries. Note that the inclusion $v_\ell T\hookrightarrow T$ also induces a homeomorphism between $\partial (v_\ell T)$ and the clopen subset $U_{v_\ell}$ of $\partial T$ that consists of all paths passing through the vertex $v_\ell$. If $v_\ell = e_\ell$, then $U_{v_\ell} = U_\ell$  for the set $U_\ell$ defined in \eqref{eq:uk}.

Denote by $\G_{v_\ell}$ the subgroup of elements $g \in \G$ which fix $v_\ell$ and so preserve $U_{v_\ell}$. For $g \in \G_{v_\ell}$ denote by $g_\ell = g|U_{v_\ell}$ the restriction, and consider the pullback $p_\ell^* \, g_\ell$ to $\partial T$. We refer to \cite{Grigorchuk2011} for a precise definition of a self-replicating group, but it implies that for any $v_\ell \in V_\ell$ and any $\ell \geq 0$, the morphisms
\begin{equation}\label{eq-replicating}
\widetilde{p}_\ell \colon  \G_{v_\ell} \to Aut(T) : g \mapsto p_\ell^* \, g_\ell = p_\ell^* (g|U_{v_\ell}) \ 
\end{equation}
have image in $\G \subset Aut(T)$ and are surjective onto $\G$. We are now in the position to establish the connection between self-replicating groups and renormalizable actions.

\begin{prop}\label{prop-selfrepl}
Let $\G$ be a self-replicating group acting on the boundary $\partial T$ of a regular tree $T$. Then the action $(\partial T,G)$ is renormalizable in the sense of Definition \ref{defn-renormalizable}.
\end{prop}

\proof In the notation of above, fix the vertex  $v_\ell = e \cdot \G_\ell$, so $U_{v_\ell} = U_\ell$ and $\G_{v_\ell} = \G_\ell$. The group $H_{U_\ell} \subset \Homeo(U_\ell) $ is  a quotient  group of $\G_\ell$. The pullback map $p_\ell^* \,\colon H_{U_\ell} \to \G$ is clearly an injective group homomorphism, as $H_{U_\ell}$ is a group of homeomorphisms of a subtree.
The map  $p_\ell^*$  is surjective by the definition of a self-replicating group and \eqref{eq-replicating}. 
 It follows that the action of $\G$ is renormalizable with maps $\lambda = p_\ell:\partial T \to \partial (v_\ell T)$ and $\Theta=(p_\ell^*)^{-1}: \G \to H_{U_\ell}$. %However, there need not be an injective map of $\G$ onto $\G_\ell$, since elements in $H_{U_\ell}$ can be extended to $\partial T$ in multiple ways if  $\Phi_{U_\ell}: \G_\ell \to H_{U_\ell}$ has a non-trivial kernel. 
 \endproof 

An example of a group whose action is renormalizable and not quasi-analytic is the Grigorchuk group, which is known to be self-replicating \cite{Grigorchuk2011}. We refer to \cite{Grigorchuk2011} for other examples of self-replicating groups, acting on trees, and to \cite{BGN2001} for the overview of the relation between the notions of self-similar groups and other notions of renormalizability, for instance that of tilings. 

Grigorchuk \cite[Proposition 11.6]{Grigorchuk2011} showed that a countable self-replicating group $\G$ which acts freely on the boundary $\partial T$ of a tree $T$ is scale-invariant, see the Introduction for the definitions. Nekrashevych and Pete \cite{NekkyPete2011} provided examples of finitely generated scale-invariant groups that are not strongly scale-invariant. Our results can be used to strengthen Grigorchuk's result in \cite[Proposition 11.6]{Grigorchuk2011} for free actions of finitely generated self-replicating groups to show that they are strongly scale-invariant.

\begin{prop}
Let $\G$ be a finitely generated self-replicating group, and suppose the action of $\G$ on the boundary $\partial T$ of a regular tree $T$ is free. Then $\G$ is strongly scale-invariant.
\end{prop}

\proof By the argument above the action of a self-replicating group is renormalizable, and if it is free, then it is quasi-analytic. Then by Proposition \ref{prop-qar2} $\G$ is strongly scale-invariant, that is, there is a renormalization $\varphi: \G \to \G$, and the group chain $\{\varphi^\ell(\G) \mid \ell \geq 0\}$ has trivial intersection.
\endproof

It is natural to ask if the converse of Proposition \ref{prop-selfrepl} holds, that is, if an effective renormalizable action is always that of a self-replicating group of homeomorphisms. In Definition \ref{defn-renormalizable} of a renormalizable action the map $\Theta: \G \to H_U$ is allowed to be any group isomorphism, introducing a possibility of a `twist'.  We leave the question whether the converse holds as an open problem.

\begin{prob}
Let $(\fX,\G,\Phi)$ be an effective renormalizable Cantor action,  and suppose $(\fX,\G,\Phi)$ is conjugate to an action of $\G$ on $(\partial T, \G)$, where $\partial T$ is the boundary of a regular tree $T$. Prove that $\G$ is a self-replicating group, or find a counterexample to this statement.
\end{prob}
  
\section{Continuous orbit equivalence}

We next give the proofs of Theorems~\ref{thm-renormal} and   \ref{thm-Dcoe}, which consider the properties of renormalizable actions which are preserved by continuous orbit equivalence.

\subsection{Proof of  Theorem~\ref{thm-renormal}}  

Let $(\fX, \G, \Phi)$ and $(\fX', \G', \Phi')$ be minimal equicontinuous Cantor actions which are continuously orbit equivalent,  and assume that     $(\fX, \G, \Phi)$ is   renormalizable and locally quasi-analytic. We claim that  $(\fX', \G', \Phi')$ is  virtually renormalizable.

First note that by Proposition~\ref{prop-qar}, the action $(\fX, \G, \Phi)$ is  quasi-analytic, and by 
Proposition~\ref{prop-qar2}, there  exists a proper self-embedding $\vp \colon \G \to \G$ such that the action $(\fX, \G, \Phi)$ is isomorphic to the  action $(X_{\vp}, \G, \Phi_{\vp})$.  
Thus,  the Cantor actions $(X_{\vp}, \G, \Phi_{\vp})$ and $(\fX', \G', \Phi')$ are continuously orbit equivalent, where $(X_{\vp}, \G, \Phi_{\vp})$ is quasi-analytic by Theorem~\ref{thm-main1}, and stable by Corollary~\ref{cor-stable}. 
Then  Theorem~6.9 of \cite{HL2019a}   implies that $(\fX', \G', \Phi')$ is   locally quasi-analytic.

The hypotheses of  Theorem~1.5 in \cite{HL2018b} are then satisfied, 
so   that $(X_{\vp}, \G, \Phi_{\vp})$ is return equivalent to $(\fX', \G', \Phi')$.
Thus, there exists adapted sets $V   \subset X_{\vp}$ for the action $(X_{\vp}, \G, \Phi_{\vp})$ and $V' \subset \fX'$ 
for the action $(\fX', \G', \Phi')$, so that the restricted actions $(V, H_V, \Phi_V)$ and $(V', H'_{V'}, \Phi'_{V'})$ are isomorphic,  where  $H_V = \Phi_V(\G_V) \subset \Homeo(V)$ and  $H'_{V'} = \Phi'_{V'}(\G'_{V'}) \subset \Homeo(V')$.

 Let $x_{\vp} \in X_{\vp}$ denote the fixed-point for the   contraction $\lambda_{\vp} \colon X_{\vp} \to X_{\vp}$.
The action  $(X_\vp, \G , \Phi_\vp)$ is minimal, so by conjugating by an element of  $\G$, we can assume that $x_{\vp} \in V$.

Let $h \colon V \to V'$ be a homeomorphism  and  $\Theta \colon H_V \to H'_{V'}$ a group isomorphism  which realizes the isomorphism between $(V, H_V, \Phi_V)$ and $(V', H'_{V'}, \Phi'_{V'})$ as in Definition~\ref{def-isomorphism}.

 For the action $(X_{\vp}, \G, \Phi_{\vp})$, we have an adapted neighborhood basis $\{U_{\ell} = \lambda^{\ell}(X_\vp) \mid \ell \geq 0\}$ and a group chain $\cG_{\vp} = \{\G_{\ell} = \vp^{\ell}(\G) \mid \ell \geq 0\}$  as before.

Choose $\ell_0 > 0$ sufficiently large so that  $U_{\ell_0} \subset V$  and $h(U_{\ell_0}) \subset V'$. 
Then  set $W = U_{\ell_0}$.
Note that $\lambda_{\vp}(U_{\ell}) = U_{\ell+1}$ for all $\ell \geq 0$, so $W_1 = \lambda_{\vp}(W)   \subset W$. 
Set $W' = h(W) \subset V'$ and $W_1' = h(W_1) \subset W'$.
Then the restriction of $\vp$ to $\G_W = \G_{\ell_0}$ yields  a proper self-embedding   $\vp_W \colon \G_W \to \G_W$. %  where we recall that  $\G_W \cong \G$.

  Since the action $(X_{\vp}, \G, \Phi_{\vp})$ is quasi-analytic, the   map $\Phi_W \colon \G_W \to H_W$ is an isomorphism. Thus, $\vp_W$ induces a proper self-embedding 
  $\whvp_W \colon H_W \to H_W$. Then set $H_{\ell} = \whvp_W^{\ell}(H_V)$ for all $\ell \geq 0$. 
  It   then follows from the constructions  that the Cantor action $(W, H_W, \Phi_W)$ is isomorphic with the Cantor action associated to $\whvp_W \colon H_W \to H_W$.

Finally, the isomorphism between $(V, H_V, \Phi_V)$ and $(V', H'_{V'}, \Phi'_{V'})$ restricts to an isomorphism between 
$(W, H_W, \Phi_W)$ and $(W', H'_{W'}, \Phi'_{W'})$ which then defines a self-embedding of $H'_{W'}$. 
Thus, the Cantor action $(\fX', \G', \Phi')$ is virtually renormalizable. 
This completes the proof of Theorem~\ref{thm-renormal}.

\subsection{Proof of  Theorem~\ref{thm-Dcoe}}
  Let $(X_{\vp}, \G, \Phi_{\vp})$ and   $(X'_{\vp'}, \G', \Phi'_{\vp'})$ be Cantor actions associated to renormalizations 
$\vp \colon \G \to \G$ and $\vp' \colon \G' \to \G'$, respectively. 
Assume that  $(X_{\vp}, \G, \Phi_{\vp})$ and   $(X'_{\vp'}, \G', \Phi'_{\vp'})$  are continuously orbit equivalent.  We must show that the discriminant groups $\cD_{\vp}$ and $\cD'_{\vp'}$ for these actions are isomorphic.  

By  Corollary~\ref{cor-stable}, the   actions    $(X_{\vp}, \G, \Phi_{\vp})$ and   $(X'_{\vp'}, \G', \Phi'_{\vp'})$ are quasi-analytic and stable. Then 
  Theorem~1.5 in  \cite{HL2018b}   implies that the actions  $(\fX, \G, \Phi)$ and   $(X'_{\vp'}, \G', \Phi'_{\vp'})$   are   return equivalent.

Thus, there exist  adapted sets $V   \subset \fX$ for the action $(\fX, \G, \Phi)$ and $V' \subset \fX'$ 
for the action $(\fX', \G', \Phi')$ so that the restricted actions $(V, H_V, \Phi_V)$ and $(V', H'_{V'}, \Phi'_{V'})$ are isomorphic, where recall that $H_V = \Phi_V(\G_V) \subset \Homeo(V)$ and  $H'_{V'} = \Phi'_{V'}(\G'_{V'}) \subset \Homeo(V')$. As the actions are quasi-analytic, the maps $\Phi_V$ and $\Phi'_{V'}$ are monomorphisms, hence are isomorphisms. Thus, the actions $(V, \G_V, \Phi_V)$ and $(V', \G'_{V'}, \Phi'_{V'})$ are isomorphic, induced by a homeomorphism $h \colon V \to V'$.

Let $\cD_V$ denote the discriminant group for the restricted action $(V, \G_V, \Phi_V)$. Then by the arguments in \cite[Section~4]{HL2018a}, 
there is a surjective map $\rho_{\fX,V} \colon \cD_{\vp} \to \cD_V$ which is an isomorphism when the profinite action $\whPhi_{\vp} \colon G \times X_{\vp} \to X_{\vp}$ is quasi-analytic. 
Likewise,  for the discriminant $\cD'_{V'}$ of the action $(V', \G'_{V'}, \Phi'_{V'})$, there is an isomorphism  $\rho_{\fX',V'} \colon \cD_{\vp'} \to \cD'_{V'}$.

The isomorphism class of the discriminant group is an invariant for isomorphism of Cantor actions, so we conclude 
$\cD_{\vp} \cong \cD_{V} \cong \cD'_{V'} \cong \cD'_{\vp'}$ as claimed.
  This completes the proof of Theorem~\ref{thm-Dcoe}.

\section{Applications and Examples}~\label{sec-examples}

In this section, we discuss some of the applications of the results of this paper, then give a selection of examples to illustrate these results.  
  
  For a compact manifold $M$ without boundary, an expansive diffeomorphism $\phi \colon M \to M$ gives rise to a renormalization $\vp \colon \G \to \G$ of the fundamental group  $\G = \pi_1(M,x)$. In this case, Shub showed in \cite{Shub1970} that the universal covering of $M$ has polynomial growth type, and hence by Gromov \cite{Gromov1981} the group $\G$ has a finite-index  nilpotent subgroup. There are a variety of  constructions of expansive diffeomorphisms   on nilmanifolds, and 
the invariants associated to the renormalization $\vp$ of $\G$ are then invariants of the expansive map $\phi$. 

   The  construction  of generalized Hirsch foliations   in \cite{BHS2006,Hirsch1975} is based on choosing a renormalization $\vp \colon \G \to \G$ of the fundamental group of a compact manifold $M$. Thus, invariants of the renormalization yield invariants for this   genre of foliated manifolds.

   The classification of $M$-like   laminations, where $M$ is a fixed compact manifold, is reduced to   the classification of renormalizations in the work \cite{CHL2018b}. 
 
These applications are all based on the   constructions of renormalizations for groups with the non-co-Hopfian property. Many finitely generated nilpotent groups are renormalizable, as shown for example in \cite{Belegradek2003,Cornulier2016,DL2003a,DL2003b,DD2016,LL2002}. There is also  a variety of  examples of renormalizable groups  which are not nilpotent, as   described for example in \cite{DP2003,ER2005,GW1992,GW1994,GLiW1994,NekkyPete2011,OP1998,WvL2018b}.
While these works show the existence of a proper self-embedding for a particular class of groups, they do not calculate the groups $\cD_{\vp}$ and $\whN_{\vp}$ which are associated to an embedding $\vp$ by Theorem~\ref{thm-main2}.
In the following, we make these calculations for a selected set of examples of renormalizable groups.

In Section~\ref{subsec-infinite} we give a basic example of a renormalizable Cantor action, where the group $\G$ is not finitely generated, and $\G$ is not virtually nilpotent.

In Section~\ref{subsec-dihedral} we give a basic example of the cross-product construction of renormalizable groups, for which the discriminant is a  non-trivial finite group.

 In Section~\ref{subsec-nilpotent} we calculate  the discriminant  $\cD_{\vp}$ and the induced map $\whvp \colon \cD_{\vp} \to \cD_{\vp}$  for an ``untwisted'' embedding $\vp \colon \cH \to \cH$ of the Heisenberg group $\cH$. 
    
In Section~\ref{subsec-padic} we give an example of a renormalizable group that   arises in the study of arboreal representations of absolute Galois groups of number fields.  

      \subsection{Infinitely generated actions}\label{subsec-infinite}

The assumption that $\G$ is finitely generated is essential for the conclusion of Theorem~\ref{thm-main4}, as shown by the following example. Let $F$ be a finite nonabelian simple group and set $\ds \G:= \oplus_{i=0}^{\infty} F$. Then $\G$ is a countable group, but not finitely generated.

Observe that $\G$ admits a renormalization, given by   by the shift map, $\vp(f_0, f_1, \ldots) = (e, f_0, f_1, \ldots)$, for $e \in F$ the identity element. The associated Cantor space is the infinite product $X_{\vp} = \prod_0^{\infty} F$.
 The action $\Phi_{\vp}$ of $\G$ on $X_{\vp}$ is free, so  the discriminant $\cD_{\vp}$ is trivial in this case. However, $\G$ is not virtually nilpotent.

       \subsection{Multihedral groups}\label{subsec-dihedral}
 
This is an elementary example of a  group $\G$ with self-embedding $\vp$ and non-trivial finite discriminant group $\cD_{\vp} \subset \G$.

Let $\Lambda = \mZ^k$ be the free abelian group on $k$ generators. Let $H \subset \Perm(k)$ be a non-trivial subgroup of  the  finite  symmetric group $\Perm(k)$ on $k$ symbols, and assume that $H$ is a simple group. Let $\Perm(k) \subset \text{GL}(k,\mathbb{Z})$ be the standard embedding permuting the coordinates.

Let $\G = \mZ^k \rtimes H$ be the semi-direct product of these groups. Fix  $m > 1$, then define $\vp \colon \G \to \G$ to be multiplication by $m$ on the $\mZ^k$ factor. That is, for $(\vec{v}, g) \in \G$ set $\vp(\vec{v}, g) = (m \cdot \vec{v}, g)$. Then 
\begin{eqnarray}
\G_{\ell} & = &  \{ (m^{\ell} \cdot \vec{v}, g) \mid \vec{v} \in \mZ^k ~ , ~ g \in H \} \   = m^\ell \mZ^k \rtimes H   \label{eq-nilchain} \\ 
K(\cG_{\vp})  & = & \{ (0, 0, g) \mid   g \in H \} \   \cong   H \ .
\end{eqnarray}
where  $\cG_{\vp} = \{\G_{\ell} \mid \ell \geq 0\}$. Then  we have 
  $X_{\vp} \cong \widehat{\mZ}_{m}^k$,  the product of $k$-copies of the inverse limit space ${\ds \widehat{\mZ}_m = \lim_{\longleftarrow}\{\mZ/m^{k+1}\mZ  \to \mZ/m^k \mZ, \, k \geq 0 \}}$.
  The subgroup  $H$ acts on  $X_{\vp}$ by permutations of the coordinates, so the adjoint action on $X_{\vp}$ of a non-identity element $g \in H$    is   a non-trivial permutation of the coordinate axes, hence is non-trivial. Thus, the normal core $C(\cG_{\vp}) \subset K(\cG_{\vp})$ is trivial, and we have $K(\cG_{\vp}) \subset \cD_{\vp}$. 
    Thus, a  calculation shows that the normal core $C_{\ell} \subset \G_{\ell}$ is the subgroup of \eqref{eq-nilchain} where $g = e \in H$ is the identity, so $\G_{\ell}/C_{\ell} \cong H$ for all $\ell > 0$. Thus, $\cD_{\vp} \cong H$.
Also, the subgroup $\whN_{\vp}$ is the   product of $k$ copies of $\widehat{\mZ}_{m}$, or the $m$-adic $k$-torus.

 Observe that the map $\vp$ restricts to the identity on the subgroup $H$, while $\vp$ acts as  multiplication by $m$ on the normal subgroup $\mZ^k$. Thus, $\whvp \colon \cD_{\vp} \to \cD_{\vp}$ in Theorem~\ref{thm-main2}, item 3,  is the identity map, and $\whvp \colon \whN_{\vp} \to \whN_{\vp}$  in Theorem~\ref{thm-main2}, item 2, is induced by coordinate-wise multiplication by $m$ on $\mZ^k$.

 \subsection{Nilpotent endomorphisms}\label{subsec-nilpotent}

The $3$-dimensional Heisenberg group $\cH$ is the simplest non-abelian nilpotent group, and we give a self-embedding  for which  $\cD_{\vp}$ is  an infinite profinite group. 

A  general construction   for self-embeddings of $2$-step nilpotent groups is given by Lee and Lee in  \cite{LL2002}, of which this example is a special case. More generally, group chains in $\cH$ were studied in detail by 
 Lightwood,  {\c{S}}ahin and  Ugarcovici in \cite{LSU2014}, where they give a complete description   for the subgroups of $\cH$ and a characterization of which subgroups are normal.   This work also gives   a   discussion of twisted and untwisted subgroups of the Heisenberg group, which can be used to construct further examples of renormalizations.

Group chains in $\cH$ whose   discriminant invariant is an infinite group were first constructed by Dyer in her thesis \cite{Dyer2015}, and also described in   \cite[Example~8.1]{DHL2016a}.  In the following, we construct such a group chain realized via a self-embedding of $\cH$.

Let $\cH$ be  represented as $(\mZ^3, *)$ with the group operation $*$, so  for $x, u,y, v,z, w  \in \mZ$ we have, 
\begin{equation}\label{eq-Hrules}
(x,y,z)*(u,v,w)=(x+u,y+v,z+w+xv) \quad , \quad (x,y,z)^{-1} = (-x, -y, -z +xy) \ .
\end{equation}
This is equivalent to the upper triangular representation in ${\rm GL(\mZ^3)}$.  In particular, we have 
\begin{equation}\label{eq-Hcomm}
(x,y,z) * (u,v,w)*(x,y,z)^{-1}=(u,v,w + xv -yu) \ .
\end{equation}

For integers $p,q > 0$ define $\vp \colon \cH \to \cH$ by a self-embedding by 
$\vp(x,y,z) = (px, qy, pqz)$. Then 
$$\cH_{\ell} = \vp^{\ell}(\cH) = \{(p^{\ell} x, q^{\ell}y, (pq)^{\ell}z) \mid x,y,z \in \mZ\} \quad, \quad \bigcap_{\ell > 0} \ \cH_{\ell} = \{e\} \ .$$
 Now assume that $p,q > 1$ are distinct prime numbers.
Formula \eqref{eq-Hcomm} implies that the normal core for $\cH_{\ell}$ is given by
$$C_{\ell} = {\rm core}(\cH_{\ell})  = \{((pq)^{\ell} x, (pq)^{\ell}y, (pq)^{\ell}z) \mid x,y,z \in \mZ\} \ .$$
Thus,  the finite group 
\begin{equation}\label{eq-Qell}
Q_{\ell} = \cH/C_{\ell} = \{( x,  y,  z) \mid x,y,z \in \mZ/(pq)^{\ell}\mZ \} \ .
\end{equation}
The profinite group $\widehat{\cH}_{\infty}$ is the inverse limit of the quotient groups $Q_{\ell}$ so we have
$$\widehat{\cH}_{\infty} =   \{(x,y,z) \mid x,y,z \in \widehat{\mZ}_{pq} \}$$
 with multiplication  on each finite quotient induced given by the formula \eqref{eq-Hcomm}.
To identify the discriminant subgroup $\cD_{\infty}$ first note   
\begin{eqnarray}
\cH_{\ell}/C_{\ell} & = & \{(p^{\ell} x, q^{\ell}y, 0) \mid x \ \in \mZ/q^{\ell}\mZ, \ y \in \mZ/p^{\ell}\mZ \} \ \subset \ Q_{\ell} \ ,  \label{eq-Qell1}\\
\cH_{\ell+1}/C_{\ell+1} & = & \{(p^{\ell +1} x, q^{\ell +1}y, 0) \mid x \ \in \mZ/q^{\ell+1}\mZ, \ y \in \mZ/p^{\ell +1}\mZ \} \   \ . \label{eq-Qell=2}
\end{eqnarray}
 
 The bonding map  $\ds q_{\ell +1} \colon \cH_{\ell +1}/C_{\ell+1} \to \cH_{\ell}/C_{\ell}$  from the definition  \eqref{eq-discformula} for $\cD_{\infty}$ is induced from the inclusion $\cH_{\ell +1} \subset \cH_{\ell}$ modulo quotient by 
$$\cH_{\ell +1} \cap C_{\ell} =    \{(p^{\ell+1} q^{\ell} x, p^{\ell} q^{\ell+1}y, (pq)^{\ell +1}z) \mid x,y,z \in \mZ\} \  .$$
 Thus, in terms of the coordinates $x,y$ in \eqref{eq-Qell=2} the bonding map is given by 
$$q_{\ell +1}(x,y,0) = (x \ {\rm mod} \ q^{\ell}\mZ, y \ {\rm mod} \ p^{\ell}\mZ,0) \ .$$
It then follows  by formula \eqref{eq-discformula} that 
\begin{equation}\label{eq-pqdisc}
\cD_{\vp} \cong  \cD_{\infty} = \varprojlim \ \{q_{\ell +1} \colon \cH_{\ell +1}/C_{\ell+1} \to \cH_{\ell}/C_{\ell} \mid \ell \geq 0\} \cong \widehat{\mZ}_{q} \times \widehat{\mZ}_{p} \ .
\end{equation}
The induced map $\whvp \colon \cD_{\vp} \to \cD_{\vp}$ is given   by multiplication by $p$ on $\widehat{\mZ}_{q}$ in the first $x$-coordinate, and multiplication by $q$ on $\widehat{\mZ}_{p}$ in the second $y$-coordinate, so that  $\whvp$ acts as an isomorphism on $\cD_{\vp}$, as asserted in Theorem~\ref{thm-main2}.

Finally, consider the subgroup of $Q_{\ell}$ in \eqref{eq-Qell} which is complementary to the subgroup $\cH_{\ell}/C_{\ell}$,
\begin{equation}
N_{\ell} = \{( q^{\ell} x, p^{\ell} y,  z) \mid  x \ \in \mZ/p^{\ell}\mZ, \ y \in \mZ/q^{\ell}\mZ, z \in \mZ/(pq)^{\ell}\mZ \} \ \subset \ Q_{\ell} \ .
\end{equation}
The   map $\vp$ induces a map on $N_{\ell}$ given  by multiplication   by $p$   in the first $x$-coordinate, and multiplication by $q$   in the second $y$-coordinate, so the action is nilpotent   on $N_{\ell}$. 
 The inverse limit of the groups $N_{\ell}$ is a subgroup of $\widehat{\cH}_{\infty}$ identified with  
$$\whN_{\vp} \cong \widehat{\cH}_{\infty}/\cD_{\infty} \cong  \{(x,y,z) \mid x \in \widehat{\mZ}_{p} \ , \ y \in \widehat{\mZ}_{q} \ ,  \ z \in \widehat{\mZ}_{pq} \}  \ , $$
and is a pro-nilpotent group as it has the finite nilpotent groups $N_{\ell}$ as quotients. Moreover,  the induced map   $\whvp \colon \whN_{\vp} \to \whN_{\vp}$ is a contraction, as asserted in Theorem~\ref{thm-main2}.
 
  Note that if we take $p=q$ in the above calculations, so $\vp \colon \cH \to \cH$ is the ``diagonal expansion'' by $p$ on the abelian factor $\mZ^2$, then $\cH_{2\ell} \subset C_{\ell}$. So while each quotient $\cH_{2\ell}/C_{2\ell}$ is non-trivial, its image under the composition of bonding maps in 
   \eqref{eq-discformula} vanishes in $\cH_{\ell}/C_{\ell}$, hence $\cD_{\vp}$ is the trivial group in the inverse limit. Correspondingly, the inverse limit space  $X_{\vp}$ has a well-defined group structure.

 \subsection{Semi-direct product of dyadic integers with its   group of units}\label{subsec-padic}
 This example  can be viewed as a more sophisticated version of Example~\ref{subsec-dihedral}.  It arises, in particular,  as the profinite arithmetic iterated monodromy group associated to  a certain post-critically finite quadratic polynomial, as discussed in \cite{Lukina2018b}.  
   
  Let $\widehat{\Gamma} = \widehat{\mZ}_2 \rtimes \widehat{\mZ}_2^\times$, where $\widehat{\mZ}_2$ is the dyadic integers, and $\widehat{\mZ}^\times_2$ is the multiplicative group of dyadic integers. Denote by $a$ the topological generator of the abelian group $\widehat{\mZ}_2$, that is, $a$ is identified with $([1]) \in \widehat{\mZ}_2$, where $[1]$ is the equivalence class of $1$ in $\mZ/2^n\mZ$, $n \geq 1$.

Recall that $\widehat{\mZ}^\times_2$ is the automorphism group of $\widehat{\mZ}_2$. 
The multiplicative units in the $2$-adic integers can be computed by computing the units in $\mZ/2^n \mZ$ for any $n$,   and taking the  inverse limit (see \cite[Theorem 4.4.7]{RZ2000}) so we have $\widehat{\mZ}^\times_2 \cong \mZ/2\mZ \times \widehat{\mZ}_2$. Here,  $\mZ/2\mZ$  is generated by $([-1]) \in \widehat{\mZ}_2^\times$, where $[-1]$ denotes the equivalence class of $-1$ in $\mZ/2^n\mZ$ for $n \geq 1$, and the the second factor is generated by $([5]) \in \widehat{\mZ}_2^\times$, where $[5]$ is the equivalence class of $5$ in $\mZ/2^n\mZ$ for $n \geq 1$. Denote these generators by $b$ and $c$ respectively. Then let 
  \begin{align}\label{eq-garith} \Gamma \cong \langle a,b,c \mid b^2 = 1, \, bab^{-1} = a^{-1}, \, cac^{-1} = a^5, \, bcb^{-1}c^{-1} = 1\rangle, \end{align}
where $b$ and $c$ commute since they are generators of different factors in a product space.

 Define a self-embedding $\vp \colon  \G \to \G$ by setting $\vp(a) = a^2$, $\vp(b) = b$ and $\vp(c) = c$.
That is, we have
  $$ \Gamma_1 = \vp(\Gamma) \cong \langle a^2,b,c \mid b^2 = 1, \, ba^2b^{-1} = a^{-2}, \, ca^2c^{-1} = (a^2)^5, \, bcb^{-1}c^{-1} = 1\rangle,$$
and so we obtain a group chain $\ds \Gamma_{\ell} = \langle a^{2^{\ell}},b,c \rangle, \, \ell \geq 1$. 
The discriminant group of the action defined by this group chain was computed in \cite[Section 7]{Lukina2018b}. In particular, computing the normal cores of the subgroups $\Gamma_{\ell}$ we obtain 
 $\ds C_{\ell} = \langle a^{2^{\ell}} ,  c^{2^{\ell}-2}\rangle \subset \Gamma_{\ell}$, 
and it follows that 
  $$\cD_\vp = \lim_{\longleftarrow}\{\Gamma_{\ell+1} / C_{\ell+1} \to \Gamma_{\ell}/C_{\ell}\}\cong \widehat{\mZ}^\times_2 \ .$$

 \section{Problems}\label{sec-problems}

The study of the   properties of the dynamical systems of the form  $(X_{\vp}, \G, \Phi_{\vp})$   suggests the following approach to the classification problem for renormalizable groups and their proper self-embeddings.

\begin{prob}\label{prob-classify}
Classify the structure of renormalizable groups $\G$ which satisfy:
\smallskip
  
 \begin{enumerate}\label{problist-classify}
\item ~ $\cD_{\vp}$ is the trivial group; 
\item ~ $\cD_{\vp}$ is   a finite group; 
\item ~ $\cD_{\vp}$ is   an infinite profinite group. 
\end{enumerate}
\end{prob}

Case (1)   is discussed further in Section~\ref{subsec-trivial} below.
There are numerous and varied constructions of examples of     case (2),  where  $\cD_{\vp}$ is a finite group. See Section~\ref{subsec-dihedral}   for some typical examples. 

The most interesting problems arise   for   case (3),  where  $\cD_{\vp}$ is an infinite profinite group.    Corollary~\ref{cor-stable} implies that all of the direct limit group  invariants for Cantor actions  defined in \cite{HL2019a} are bounded for these examples.   Thus, 
the problem is to refine the invariants constructed from   the adjoint action of $\cD_{\vp}$ on the pro-nilpotent normal subgroup  $\whN_{\vp} \subset \whGamma_{\vp}$ to distinguish these various examples. Note that if the group chain $\cG_{\vp}$ has trivial intersection,   then the intersection $\cD_{\vp} \cap \G$ is trivial, so the invariants constructed using the adjoint action of $\cD_{\vp}$ are only ``seen'' when considering the action of  the profinite group $\whGamma_{\vp}$.

 \subsection{Renormalizable nilpotent groups} \label{subsec-trivial}
 
 Suppose that $\G$ admits a renormalization $\vp \colon \G \to \G$, such that each of the subgroups $\G_{\ell} = \vp^{\ell}(\G)$ is a \emph{normal} subgroup of $\G$. Then the third author showed in the work \cite{WvL2018a} that the quotient $\G/C(\cG_{\vp})$ must be free abelian. In particular, if the group chain $\cG_{\vp} = \{\G_{\ell} \mid \ell \geq 0\}$ has trivial intersection, then $\G$ is free abelian. 
 Theorem~\ref{thm-main4} is a more general form of this result, where the assumption that $\cG_{\vp} $ has finite discriminant implies that $\G$ is virtually nilpotent. 
 
  The remarks at the end of Section~\ref{subsec-nilpotent} show that $\cD_{\vp}$ is trivial when $p=q$ for the construction in  Section~\ref{subsec-nilpotent}. In fact, these remarks apply in general   to the diagonal action on the   nilpotent subgroup of upper triangular integer matrices, where $\vp$ is given by multiplication by a constant factor $p$ on the super-diagonal entries; that is, those directly above the diagonal.   
 This suggests that   the non-triviality of the discriminant  invariant $\cD_{\vp}$ for an endomorphism of a nilpotent group is a measure of the ``asymmetry''   of the embedding $\vp$.    It is an interesting problem to make this statement more precise for the general nilpotent group.

\begin{prob}\label{prob-nilinvariants}
Let $\G$ be a finitely generated torsion free nilpotent  group,  and $\vp$ a renormalization such that $\cG_{\vp} = \{\G_{\ell} \mid \ell \geq 0\}$ has trivial intersection. Develop    the relationship between the properties of  the  discriminant group $\cD_{\vp}$, the embedding $\vp$, and the nilpotent structure theory of $\G$,  as developed for example in   \cite{Cornulier2016,DD2016}.
\end{prob}

 \subsection{Algebraic invariants}\label{subsec-algebraic}

The reduced group $C^*$-algebra  $C_r^*(X_{\vp}, \G, \Phi_{\vp})$ obtained from the group action $(X_{\vp}, \G, \Phi)$ is a source of invariants for the group $\G$ and the embedding $\vp$. In the case when $\G = \mZ^n$ is free abelian, the work \cite{GPS2019} shows that the ordered K-theory of this $C^*$-algebra is a complete invariant of the action. It is natural to ask whether    similar results are possible in more generality:
 \begin{problem}\label{prob-nilpotentinv}
Let $\G$ be a finitely generated nilpotent group, and $\vp$ a renormalization of $\G$. What information about the nilpotent structure constants of $\G$ and the embedding $\vp$ is determined by the K-theory groups $K_*(C_r^*(X_{\vp}, \G, \Phi_{\vp}))$?
\end{problem}

Note that by Theorem~\ref{thm-Dcoe}, the isomorphism class of the discriminant group $\cD_{\vp}$ is an invariant of the continuous orbit equivalence class of the Cantor action $(X_{\vp}, \G, \Phi_{\vp})$, and the isomorphism class of $C_r^*(X_{\vp}, \G, \Phi_{\vp})$ is also   invariant. It   seems natural that these two invariants should be closely related.

 \begin{problem}\label{prob-nilpotentinv2}
Let $\G$ be a renormalizable group.
How does the algebraic structure of $C_r^*(X_{\vp}, \G, \Phi_{\vp})$  reflect the properties of the profinite group $\cD_{\vp}$?
\end{problem}

Theorem~\ref{thm-main2} shows that the profinite group $\whGamma_{\vp}$ is a semi-direct product with $\cD_{\vp}$ as a factor. One approach to Problem~\ref{prob-nilpotentinv2} would be to relate the decomposition $\whGamma_{\vp} \cong \whN_{\vp} \rtimes \cD_{\vp}$ in Theorem~\ref{thm-main2} to the algebraic structure of $C_r^*(X_{\vp}, \G, \Phi_{\vp})$.

 \subsection{Realization}\label{subsec-scale}

Given any pro-finite group  $\cD$ which is topologically countably generated,   it was shown in \cite{HL2018a,HL2019a}, using the Lenstra method,  that there exists a finitely generated group $\G$ and Cantor action $(\fX, \G, \Phi)$ whose discriminant is isomorphic to $\cD$. 
 
\begin{problem}\label{prob-siscale}
  What  profinite  groups can be realized  as the discriminant  of a Cantor action  associated to a  renormalizable group $\G$?   
\end{problem}

 \subsection{Renormalizable Cantor actions}\label{subsec-virtual}

 Let  $(\fX,\Gamma,\Phi)$ be a minimal  equicontinuous  Cantor action of wild type; that is,  the action 
 $\whPhi \colon \overline{\Phi(\G)} \to \Homeo(\fX)$    is not locally quasi-analytic. The action is said to be \emph{wild of finite type} if, in addition, 
  for some $x \in \fX$ and every clopen set with $x \in U$, we that the   kernel of the restriction 
 $\ds \ker \{ \whPhi_{U} \colon \overline{\Phi(\G)}_x \to \Homeo(U)\}$. 
   Examples of wild actions constructed by the first two authors in \cite{HL2018a} are of finite type. 
 However, the examples in \cite{HL2018a} are not renormalizable.

\begin{problem}
 Do there exist renormalizable Cantor actions which are wild of finite type?
 \end{problem}
  
\begin{problem}
Suppose that $(\fX,\Gamma,\Phi)$ is a renormalizable Cantor action which is wild. What can be said about the algebraic properties of $\G$? For example, must $\G$ have exponential growth type? What can be said    about the profinite group $\overline{\Phi(\G)} \subset \Homeo(\fX)$ for such actions?
\end{problem}

 \subsection{Representations of Galois groups}\label{subsec-galois}

  The works of the second author \cite{Lukina2018a,Lukina2018b} define  the discriminant invariants associated to arboreal representations of absolute Galois groups for number fields and function fields. Such a representation is a profinite group, obtained as the inverse limit of finite Galois groups, which act on finite extensions of the ground field, obtained by adjoining the roots of the $n$-th iteration of the same polynomial, for $n \geq 1$.

 The example given in Section~\ref{subsec-padic} is an example of an arboreal representation of an absolute Galois group,  which is isomorphic to a Cantor action  associated to a renormalization. For many polynomials the associated action is known to be not locally quasi-analytic  \cite{Lukina2018b} and, therefore, by Theorem \ref{thm-main1} it cannot be associated to a renormalization of a group. This suggests the following problem:
 \begin{problem}\label{prob-galois}
For which arboreal representations of absolute Galois groups does there exists a dense finitely generated group $\G$ and a renormalization  $\vp \colon  \G \to \G$,   such that the arboreal representation of $\G$ is return equivalent  to  a Cantor action  associated to  $(X_\vp,\G, \Phi_\vp)$?
\end{problem}

Although, as discussed above, many arboreal representations are not associated to an finite-index  embedding $\vp \colon \G \to \G$, since they are associated to a structure  built using iterations of the same polynomial, it is natural to look for a formalism similar to the non-co-Hopfian setting for the study of these  groups. This motivated the definition of renormalizable actions in Section \ref{sec-renormalization}, and suggest the following interesting problem:

\begin{problem}\label{prob-galois-renorm}
Let $(\fX, \G,\Phi)$ be an equicontinuous minimal Cantor action, and suppose that $(\fX, \G,\Phi)$ is renormalizable as in Definition \ref{defn-renormalizable}. Develop a structure theory for the group obtained as the closure of the action $(\fX,\G,\Phi)$ in $\Homeo(\fX)$,  analogous to   Theorem \ref{thm-main2}.
\end{problem}
 
%%%%%%%%%%%%%%%%%%%%%%%%%%%%%%%%%%%%%%%%%%%%%%%%%%%%%%%

\end{document}